\documentclass[trsc,nonblindrev]{informs3} 

\OneAndAHalfSpacedXI 

\usepackage[section]{placeins}  

\usepackage{comment}

\usepackage{xcolor}
\usepackage{colortbl,booktabs}
\usepackage{booktabs}
\usepackage{multirow} 
\usepackage[figuresright]{rotating} 
\usepackage{pdflscape}  

\usepackage{threeparttable} 

\usepackage{verbatim}

\usepackage{algorithm,float}  
\usepackage{algorithmic}   

\makeatletter
\newenvironment{breakablealgorithm}
{
		\begin{center}
			\refstepcounter{algorithm}
			\hrule height.8pt depth0pt \kern2pt
			\renewcommand{\caption}[2][\relax]{
				{\raggedright\textbf{\ALG@name~\thealgorithm} ##2\par}%
				\ifx\relax##1\relax 
				\addcontentsline{loa}{algorithm}{\protect\numberline{\thealgorithm}##2}%
				\else 
				\addcontentsline{loa}{algorithm}{\protect\numberline{\thealgorithm}##1}%
				\fi
				\kern2pt\hrule\kern2pt
			}
		}{
		\kern2pt\hrule\relax
	\end{center}
}
\makeatother

\usepackage{natbib}
\bibpunct[, ]{(}{)}{,}{a}{}{,}%
\def\bibfont{\small}%
\TheoremsNumberedThrough     

\EquationsNumberedThrough    

\usepackage{cleveref}
\crefname{section}{§}{§§}
\Crefname{section}{§}{§§}

\begin{document}
	
	
	\RUNAUTHOR{Huang et al.}
	
	\RUNTITLE{Integrated Air Cargo Recovery}
	
	\TITLE{A machine learning based column-and-row generation approach for integrated air cargo recovery problem}
	
	\ARTICLEAUTHORS{%
		\AUTHOR{Lei Huang, Fan Xiao, Zhe Liang}
		\AFF{School of Economics and Management, Tongji University \EMAIL{leihuang@tongji.edu.cn,fanxiao@tongji.edu.cn,liangzhe@tongji.edu.cn}, \URL{}}
	} 
	
	\ABSTRACT{%
		Freighter airlines need to recover both aircraft and cargo schedules when disruptions happen. This process is usually divided into three  sequential decisions to recovery flights, aircraft, and cargoes. This study focuses on the integrated recovery problem that makes aircraft and cargo recovery decisions simultaneously. We formulate two integrated models based on the flight connection network, one is the arc-based model, and the other is the string-based model. The arc-based model makes the flight delay decisions by duplicating flight copies, and is solved directly by commercial solvers such as Cplex. The string-based model makes the flight delay decisions in the variable generation process. The main difficulty of the string-based model is that the number of constraints grows with the newly generated flight delay decisions. Therefore, the traditional column generation method can not be applied directly. To tackle this challenge, we propose a machine learning based column-and-row generation approach. The machine learning method is used to uncover the critical delay decisions of short through connections in each column-and-row generation iteration by eliminating the poor flight delay decisions. We also propose a set of valid inequality constraints which can greatly improve the objective of LP relaxation solution and reduce the integral gap. The effectiveness and efficiency of our model is tested by simulated scenarios based on real operational data from the largest Chinese freighter airlines. The computational results show that a significant cost reduction can be achieved with the proposed string-based model in reasonable time.

	}%
	
	
	\HISTORY{}
	
	\maketitle
	
	%
	
	
		
	\section{Introduction}
	\subsection{Background}
	Air cargo has played an essential role in world trade and has an average growth of 4.1\% in volume since 1989 \citep{Boeing2020}. Although it only makes up less than 1\% proportion of global trade in terms of volume, air transport represents around 33\% of global trade in terms of value \citep{IATA2019}. Despite the fact that the COVID-19 pandemic damaged the market in 2020, a recovery is taking hold and it is forecast to grow at 4.0\% per year over the next 20 years \citep{Boeing2020}. 
	In air cargo industry,  goods are transported either in dedicated freighters, which refers to aircraft designed exclusively for cargo transport, or in the belly space of passenger aircraft \citep{Feng2015}. Freighters have clear advantages over belly space because they are separate from passenger transportation business and highly controllable. 
	In fact, freighters carry more than 50\% of the world’s air cargo traffic and generate nearly 90\% of the total air cargo industry revenue. From 2020 to 2039, the number of freighters is predicted to increase by approximately 60\%, from 2,010 to 3,260 \citep{Boeing2020}.
	
	In the real world, disruptions bring more challenges for the air operation for both passenger and cargo transport. 
	Disruptions are mainly caused by unexpected resource shortages, such as aircraft unavailability owing to mechanical issues, crew absence due to illness or disruptions of upstream schedules, as well as capacity shortage at airports and/or airspace during severe weather conditions \citep{Liang2018}.
	In the air cargo business, freighter airlines are also faced with capacity mismatch caused by demand fluctuation in the short run \citep{Delgado2020,Sandhu2006, Feng2015}. 
	When disruptions occur, flights cannot take off as scheduled, cargo routes are disrupted, or aircraft need to be rescheduled to accommodate changes in cargo demand. The Airline Operations Center (AOC) is responsible for rescheduling aircraft, flights, crews as well as cargoes/passengers to recover the airlines’ operation with minimized cost. 
	
		
	Although the recovery operations look similar, there are some non-negligible differences between the air cargo recovery problem and the passenger recovery problem (see Table \ref{Table:CargoVsPsg}). 
	The first difference is that when disruptions occur, each passenger's itinerary should be recovered as early as possible. For cargoes, however, it is sufficient as long as they can be delivered before the deadline according to the contracts.
	Second, in contrast to passenger traffic, cargo has no strong preference for a specific itinerary as long as its commitment is satisfied. In recovery process, airlines usually do not change passenger itineraries if no flights are canceled.
	On the contrary, it is much more flexible and inexpensive for freighter airlines to change the pre-scheduled cargo itineraries.
	For cargoes, the total travel time consumed and the number of connections passed are not critical.
	For example, cargo can be re-routed from a direct itinerary to an itinerary including multiple transshipment. 
	%
	Third, cargo can only be transshipped at hubs and requires additional connection time due to the specific transshipment equipment, whereas passengers can transfer at any airport, making cargo transshipment a more challenging problem than passenger transfer.
	Forth, while making flight delay decisions, passenger airlines mainly consider aircraft or crew connections but rarely consider passenger connections. On the contrary, freighter airlines are more likely to wait longer for cargoes because the revenue impact for missing cargo is much higher. 	
	
	In particular, two situations might occur in case of flight delay. In the first situation, if the two connected flights of a cargo itinerary are operated by different aircraft and the previous flight is delayed, the second flight should be delayed to catch up the cargo. Otherwise, the cargo will miss the connection. As we can see, because the cargo itinerary links two different aircraft, the delay is propagated from one aircraft to another. In the second situation, if both flights of a cargo itinerary are operated by the same aircraft, the delay is retained without affecting other aircraft. We call this type of connection as through cargo connection. Furthermore, if the through connection time is shorter than the standard cargo transshipment time between different aircraft, we call it a short through cargo connection. There are plenty of advantages to promote the through connections and short through connections. However, by considering the through connection and short through connection, delay decisions become more complicated and sophisticated to make.
	
	\begin{table}[htbp]
		\centering
		\caption{\label{Table:CargoVsPsg} Differences between Cargo Recovery and Passenger Recovery}
			\begin{tabular} 
				{l|c|c|c|c}
				\toprule
				& with Deadline
				& Change Itinerary
				& Transshipment
				& Delay for Psg/Cargo
				\\
				\midrule		
				\textbf{Passenger} 
				& No 
				& Seldom 
				& Hub/Spoke 
				& No \\
				\textbf{Cargo} 
				& Yes 
				& Yes 
				& Hub only
				& Yes \\
				\bottomrule
			\end{tabular}
	\end{table}%

	
	Compared with the planning stage, a critical constraint for airline recovery problems is the short permissible time limited to obtain solutions. 
	Freighter airlines usually utilize sequential approaches with aircraft and cargo schedule recovery rather than solving an integrated model due to the quick response requirements. In other words,  flights and aircraft are rescheduled first, and cargoes are rescheduled later according to the updated flight schedule. Because cargo rerouting strategies are not considered in the flight recovery stage,  this approach could lead to more cargo disruptions. Furthermore, the independent flight recovery failed to take into account individual emergency.
	Therefore, we suggest that an integrated model should produce a better recovery solution with a lower total recovery cost in light of these problems.
	
	This study aims to develop a model to assist freighter airlines in managing unpredictable schedule disruptions. We focus on the topic of recovery for freighters because they transport the majority of the valuable cargo. An integrated model is proposed to address aircraft and cargo recovery issues simultaneously in order to overcome the drawbacks of sequential recovery and achieve at a recovery solution with the lowest total recovery costs.
	
	\subsection{Literature Review} \label{LiteratureReview}
	
	
	There has been rich literature in airline disruption management for passenger airlines. We refer to a recent review by \cite{Su2021}. 
	Most existing literature builds the recovery model based on two types of model presentations, the time-space network and the flight connection network.
	Firstly, the time-space network (also known as the time-line or  activity-on-edge network), is a very popular approach that was first used by \cite{Yan1996} 
	to deal with the aircraft recovery problem.
	This network presentation strategy is also used in latter studies for integrated recovery problems \citep{Sinclair2014,Arkan2016,Marla2017}.
	A recent study by \cite{Huang2021} introduced a flight copy evaluation method based on the time-space network. This method can generate a limited number of flight copies to reduce recovery costs and is approved to provide promising recovery solutions in a respectable amount of time.
	Secondly, based on the flight connection network (also known as activity-on-node network), some studies utilize flight connection arcs as decision variables \citep{Arkan2017} while others use string-based variables. 
	The concept of flight string is first introduced by \cite{Barnhart1998}, it is a sequence of flights with timing decisions operated by the same aircraft. String-based models are capable of capturing network effects that individual flight decisions do not. 
	Despite the fact that the number of strings naturally expands dramatically with the number of flights, researchers have offered effective column generation (or column-and-row generation) solution approaches that can limit the number of flight strings as well as the size of the problem \citep{Petersen2012,Maher2016,Liang2018}.
	
	As to the integrated airline recovery problem, many researches have proposed models of aircraft recovery integrated with crew recovery \citep{Maher2016}, passenger recovery \citep{Arkan2016,Marla2017}, or both \citep{Petersen2012,Maher2015,Arkan2017}. In comparison to rich research on passenger airlines, literature on air cargo is relatively scarce. The integrated recovery problem of schedule, aircraft and cargoes is also merely addressed. Considering the similarity of the problems regardless of the recovery entities, we refer to the integrated aircraft recovery and passenger recovery model as a reference first, and discuss the scheduling researches for air cargo industry later.
	
	Aircraft recovery integrated with passenger recovery has received increasing attention in recent literature, and is first addressed by \cite{Bratu2006}. They proposed two optimization models to minimize jointly airline operating costs and estimated passenger delay and disruption costs, one is Disrupted Passenger Metric (DPM) only consider flight disruptions and the other is Passenger Delay Metric (PDM) with passenger re-assignment decision. The two models are solved under simulated situations with OPL Studio. 
	
	\cite{Petersen2012} studied an integrated airline recovery problem with mixed-integer for the schedule, aircraft, crew, and passenger recovery for a single-day horizon and solved the problem using Benders decomposition and column generation methods.
	\cite{Maher2015} modeled passenger recovery by prescribing alternative travel arrangements for passengers in flight cancellations circumstances and solved the integrated problem with a column-and-row generation approach.
	\cite{Arkan2016} 
	proposed an integrated aircraft and passenger recovery problem with controllable cruise speed. Cruise speed control is "a two-edged sword" since it reduces flight delay on the one hand and increases fuel cost on the other hand. 
	\cite{Arkan2017} further extended their research to include crew recovery. \cite{Marla2017} proposed a model considering flight planning with discrete flight speeds and studied the trade-off between delays and fuel burn.
	
	Researchers also presented a number of heuristic methods to obtain a fast and efficient solution. \cite{Jafari2010} introduced an integrated aircraft and passenger recovery model for recovering both aircraft and passengers simultaneously, and solved the problem using a heuristic approach. \cite{Bisaillon2011} solved the integrated problem by a large neighborhood search heuristic approach, and their work is further extended by \cite{Sinclair2016}. Other heuristic methods can refer to New Connections and Flights (NCF) \citep{Jozefowiez2013} and greedy randomized adaptive search (GRASP) \citep{Hu2016}. However, heuristic methods may result in solutions with big gap compared with the optimal solution.
	
	Existing literature in air cargo operation mainly focuses on transshipment airport selection, fleet assignment, flight routing, and cargo routing \citep{Feng2015}. The first piece of literature on planning for the air cargo industry is presented by \cite{Marsten1980}. They solved the problem of Origin-Destination (OD) pairs selection and freighter fleet assignment. Using benders decomposition, \cite{Li2006} suggested and solved the integrated problem of fleet assignment and cargo routing. \cite{Yan2006} introduced an integrated scheduling model for airport selection, fleet routing, and timetable setting in the short-term. In their study, the flight route for each airplane is generated after the fleet decision. \cite{Derigs2009} formulated an integrated model that simultaneously optimizes aircraft rotations and cargo routes and builds a solution procedure using column generation with shortest path algorithms. \cite{Derigs2013} further studied the problem with integrated fleeting and aircraft rotation. Recently, \cite{xiao2022integrated} proposed an arc-based model and a string-based model to solve the integrated aircraft and cargo routing problem and introduced the benefit of "short through cargo connections".
	
	Regarding the literature on cargo recovery, \cite{Delgado2020} and \cite{Delgado2021} addressed the problem of schedule and routing redesign of aircraft and cargo to deal with demand fluctuation circumstances in the short run. Cargo transport has higher uncertainty than passenger transport in terms of capacity availability because cargo bookings are made in very short time windows \citep{Sandhu2006}, and freight forwarders usually do not need to pay for unused capacity as well as reservation change \citep{Feng2015}. However, to the best of our knowledge, the existing studies have not addressed the recovery policies of continuous flight delays. Schedule disruptions are not addressed in the literature on cargo recovery, either.
	
	\subsection{Contribution} \label{Contribution}
	Although the integrated airline recovery problem has been extensively studied for passenger airlines, research on disruption management for freighters is scarce. As previously stated, the business of freighter airlines exhibits some differences from that of passenger airlines. Thus the existing models and solution algorithms are incapable of dealing with disruption scenarios for freighters.
	In this study, we propose an integrated model for aircraft and cargo recovery. We introduced both arc-based and string-based flight connection models and proposed an integral column-and-row with a machine learning prediction model to solve the string-based model. The main contributions of this study that distinguish it from prior studies can be summarized as follows. 
	
	First, we propose an integrated recovery solution for both aircraft and cargo under disruption circumstances instead of solving them sequentially. Cargo re-routing options are more flexible than passenger recovery, which makes the problem more complex. We also determine the delivery priority for each cargo according to the deadline in the contract. Moreover, short through connection constraints are also considered in the integrated model. Although originally scheduled short through connections might be broken when disruptions happen, freighter airlines can utilize the benefits of short through connections by generating new ones.
	
	Second, we provide an efficient column-and-row generation approaches for solving the string-based flight connection model.
	This approach solves the aircraft re-routing sub-problems as well as cargo re-routing sub-problems simultaneously, resulting in a high-quality recovery solution in a short period of time.
	Because freighter airlines are willing to delay longer for cargo connections as mentioned above, we make flight delay decisions for both aircraft and cargoes.
	The recovery algorithm addressed in this study may also shed some light on the recovery problem under passenger connection brokerage circumstances.
	
	
	Third, we integrated a machine learning algorithm to the column-and-row generation to promote the quality of  beneficial flight delay decisions and control the problem size.
	Specifically, we predict the probability of critical flight delay decisions using a decision tree algorithm, which is trained using historical data. Then we add the promising flight delays into the integrated model. This strategy reduces the number of rows added to the model, which improves the efficiency of the column-and-row generation approach.
	
	The remainder of this paper is organized as follows:  \cref{Sec:ProblemDefine} first introduces the disruptions in air cargo operation as well as the recovery decisions. We also present the problem formulation for both arc-based model and string-based model. In \cref{Sec:ColRow}, we detail the column-and-row-generation solution approach, the short through connection prediction method as well as the integral column-and-row generation with prediction algorithm. Several computational studies are presented in \cref{Sec:Experiments}. In the last section, we review our approach and present some further research directions.
	
	\section{The Air Cargo Recovery Problem(ACRP)} \label{Sec:ProblemDefine}
	\subsection{Problem Definition}

	Similar to passenger airlines, freighter airlines may encounter disruptions in operations, such as  Aircraft-On-Ground (AOG), airport/airspace closure, crew absence, and so on.
	Except for the above disruptions, freighter airlines also accept a large portion of emergency orders which is close to the flight departure. Thus they have to reschedule the original schedule to adapt to demand changes.
	Therefore, when disruptions happen, freighter airlines have to redesign an operational schedule for aircraft operation and cargo shipment schedule within a recovery horizon. This is known as the air cargo recovery problem (ACRP). 
	The problem consists of three primary issues that need to be resolved: 
	
	1) The \textit{flight recovery problem} aims to repair the disrupted original flight schedule by flight cancellation as well as flight re-timing; 
	
	2) The \textit{aircraft recovery problem} concentrates on re-rerouting each aircraft to match the restored schedule;
	
	3) The \textit{cargo recovery problem} is to re-accommodate disrupted cargoes to new itineraries that transport them to their destination before the contracted delivery time, if possible.
	 
	For ACRP in this study, we make recovery decisions for aircraft and cargoes simultaneously. We consider the following recovery policies in the integrated problem. As for flight and aircraft recovery, three recovery options are considered: flight swap, flight delay, and flight cancellation. The recovery schedule is better when it involves fewer flight cancellations, shorter flight delays, and fewer flight swaps. 
	As for cargo recovery, we assume that cargoes can not be dropped in the middle of the itineraries. That is, cargoes are either transported from the origination to the destination with the original or new generated itinerary, or they are canceled.
	The cargo recovery options include cargo re-routing, cargo delay and cargo cancellation. If the cargo is delivered to the destination after the contracted delivery time, cargo delay penalty occurs corresponding to the length of the delay. 
	Furthermore, the original cargo itinerary may become unavailable because of capacity limitation due to flight cancellation or demand fluctuation. In such circumstances, if new feasible itineraries cannot be found, the cargo need to be canceled whole or in part.
	
	When the flight and aircraft schedule change, it could cause the original cargo itineraries to be infeasible. Specifically, in cargo itineraries, a \textit{cargo connection} refers to a pair of flights that are carried sequentially. A \textit{through connection} requires the cargo staying in the same aircraft for two connected flights; hence, the cargo connection time is identical to the aircraft turn time. A through connection is called a \textit{short through connection} if the connection time is less than the minimum cargo transshipment connection time between aircraft \citep{xiao2022integrated}.
	In cargo recovery, the originally scheduled short through connections may be broken by disruptions. The cargo then must be re-assigned to new itineraries. 

	In this paper, instead of considering weight and volume, we consider the capacity as the number of unit load devices (ULDs, e.g., container or pallet). This is due to the fact that, in practice, all cargoes are pre-packed in containers and pallets with certain standards. Weight and volume constraints are no longer the bounding constraints because they are partially ensured during the container parking, loading, and balancing procedures. According to \cite{brandt2019}, most flights are not operated close to the aircraft weight capacity, and the physical volume capacity can only be filled up to 60 or 70\%. On the contrary, the constraints on the number of total containers/pallets become very critical. Thus, we utilize the number of ULDs to measure the aircraft capacity as well as cargo volume in this study. 
	
	To summarize, the ACRP seeks to obtain an optimal recovery plan with minimized overall operation costs of the recovery policies for both aircraft re-routing and cargo re-assignment. 
	We have mentioned three typical representations used in recovery problem in \cref{LiteratureReview}, namely the time-space network representation, as well as the arc-based and string-based representation on flight connection network. Since it is not easy to specify short through connections with the time-space network structure, we propose our model based on the flight connection network, with both arc-based and string-based representation.
	
	\subsection{Arc-based Model (ACRP-A)}
	
	We first introduce an arc-based model on the flight connection network, and denote it as model ACRP-A. Before we present the model, we describe the total notations in Table \ref{ConnectionNotation}.
	We first built a flight connection network for aircraft and cargoes, represented as $G$. 
	The flight delay recovery policy is implemented using flight copies. Each flight is duplicated according to a set of later departure time, represented as flight nodes in the network. 
	For each aircraft $a$,  we introduce a source node $n_a^-$ that connects to the flight departs from where the aircraft is located at the beginning of the recovery horizon, and introduce a sink node $n_a^+$ that is connected to all the flights.
	For each cargo $o$, we introduce a source node $n_o^-$ that connects to flights that depart from the cargo's origination, and the sink node $n_o^+$ that is connected to flights that arrive at the cargo's destination. Edges in the connection network, denoted as $E$,  represent the feasible flight connections between each pair of nodes. A flight connection is feasible if the minimum connection time condition is met, which depends on the type of the flight as well as the throughput of the airport.
	The set of short through connections, which is a subset of $E$, is given by $E_{sc}$ .

	The recovery task in ACRP-A is to select a set of optimal edges out of the edge set. The decision variable $u_{ij}^a$ equals to one if aircraft $a$ flows through arc $(i,j)$, at which point the flight nodes $i$  and $j$ are both selected. It determines the choice of flight delay and aircraft swap decision for the two fights.
	Similarly, the decision variable $v_{ij}^o$ states the volume of cargoes shipped from node $i$ to node $j$.
	Thus flight change compared with the original itinerary is determined. The cargo delay decision is related to arcs connected to the sink node, which is $v_{in_o^+}$.
	Additionally, the cancellation for flight $f$ is denoted as $y_f$, and the cancellation amount of cargo $o$ is denoted as $z_o$.
	In circumstances where the airport is not available, the decision variables related to the edges are equal to zero.
	
	
	\begin{table}[htbp] 
		\centering
		\small
		\caption{\label{ConnectionNotation}List of notations for model ACRP-A} 
		\begin{tabular}{p{5em}p{30em}} 
			\toprule 
			\textbf{Sets} &\\
			$A   $& Set of aircraft indexed by $a$;\\
			$F   $& Set of flights indexed by $f$;\\
			$O   $& Set of cargoes indexed by $o$;\\
			$G	$& The connection network for aircraft and cargoes;\\	
			$N$ & Set of flight copy nodes in the connection network G, indexed by $i$ and $j$;\\
			$N_f$ & Subset of nodes $N$, representing copy node set for flight $f$;\\
			$n_a^-$		& Dummy source node for aircraft $a$;\\
			$n_a^+$		& Dummy sink node for aircraft $a$;\\
			$n_o^-$		& Dummy source node for cargo $o$;\\
			$n_o^+$		& Dummy sink node for cargo $o$;\\
			$N_a^*$&$ N \cup \{n_a^-,n_a^+\}$;\\
			$N_o^*$&$ N\cup \{n_o^-,n_o^+\}$;\\
			$E$ & Set of edges in the connection network G, indexed by $(i,j)$;\\
			$E_{sc}$ & Subset of edges $E$, representing short connect edges;\\
			\midrule
			\textbf{Parameters} &\\
			$c_i^a$		& Cost of assigning flight copy $i$ to aircraft $a$;\\  
			$c_{f}  $	& Cost incurred if flight $f$ is canceled;\\
			$c_i^o $	& Cost of assigning flight copy $i$ to cargo $o$;\\  
			$c_{o}  $& Cost for every unit of cargo $o$ canceled;\\
			$d_{o}   $& Amount of cargo $o$;\\
			$ Cap_{a}    $& Capacity of aircraft $a$;\\
			\midrule
			\textbf{Variables}&\\
			$u_{ij}^a$		& $u_{ij}^a = 1$, if aircraft $a$ flies through connection $i\rightarrow j$; otherwise, $u_{ij}^a=0$;\\
			$y_f    $	& $y_f=1$ if flight $f$ is canceled; otherwise, $y_f=0$;\\
			$v_{ij}^o$		& Amount of cargo $o$ transported on connection $i\rightarrow j$;\\
			$z_{o} 	$	& Amount of cargo $o$ canceled.\\
			\bottomrule 
		\end{tabular} 
	\end{table}
	
	The full formulation of model ACRP-A is presented as follows:

	\begin{align}
		\min \ &\sum_{f \in F} c_{f} y_{f} + \sum_{a \in A}\sum_{i\in N_a^* }\sum_{j\in N_a^*}c_i^a u_{ij}^a + \sum_{o \in O} c_{o} z_{o} +  \sum_{o \in O}\sum_{i\in N_o^* }\sum_{j\in N_o^*}c_i^o v_{ij}^o  	\label{ConnObj}\\
		&\sum_{j \in N_a^*} u_{ij}^a = 1, \forall a\in A, i=n_a^-				\label{ConnAftBalance1}\\
		&\sum_{j \in N_a^*}u_{ji}^a - \sum_{j \in N_a^*} u_{ij}^a = 0, \forall a\in A,\forall i \in N		\label{ConnAftBalance2} \\
		&\sum_{i \in N_a^*} u_{ij}^a = 1, \forall a\in A, j = n_a^+							 				\label{ConnAftBalance3}\\
		&\sum_{a \in A} \sum_{i \in N_f}\sum_{j \in N_a^*} u_{ij}^a +y_{f}=1, \forall f \in F 			\label{ConnFltCover}\\
		& \sum_{j \in N_o^*} v_{ij}^o = d_o - z_o, \forall o \in O , i = n_o^-							\label{ConnCargoBalance1}\\
		& \sum_{j \in N_o^*}v_{ji}^o - \sum_{j \in N_o^*} v_{ij}^o = 0,  \forall o \in O, i \in N 			\label{ConnCargoBalance2}\\
		& \sum_{i \in N_o^*} v_{ij}^o = d_o - z_o, \forall o \in O , j = n_o^+		\label{ConnCargoBalance3}\\
		&\sum_{a \in A}\sum_{j \in N_a^*}Cap_a u_{ij}^a \ge \sum_{o \in O}\sum_{j \in N_o^*} v_{ij}^o , \forall i \in N
		\label{ConnLegCap}\\
		&\sum_{a \in A}Cap_a u_{ij}^a \ge \sum_{o \in O} v_{ij}^o , \forall (i,j) \in E_{sc}		\label{ConnSC}\\
		&u_{ij}^a \in\{0,1\} , 			\ \forall a \in A, \forall i \in N_a^*,\forall j \in N_a^*
		\label{ConnX}\\
		&y_f \in\{0,1\} ,\ \ \ \ \ \ \forall f \in F
		\label{ConnY}\\
		&z_{o}\in \mathbb{Z}, 0 \leq z_{o} \leq d_o ,	\ \forall o \in O
		\label{ConnZ}\\
		&v_{ij}^o \in \mathbb{Z} , \forall o \in O, \forall i \in N_o^*,\forall j \in N_o^* \label{ConnW}
	\end{align}
	
	The objective function (\ref{ConnObj}) minimizes the overall recovery cost, which is calculated as a summation of recovery policy costs for aircraft/flight and cargoes. Aircraft recovery costs includes cost of flight cancellation, aircraft swap and flight delay cost. Cargo recovery costs include cost of cargo cancellation, flight change and cargo delay.
	Flow balance constraints for each aircraft is defined by	constraints (\ref{ConnAftBalance1}) -- (\ref{ConnAftBalance3}). The constraints for cover of each flight is	equation (\ref{ConnFltCover}). Flights are either covered or canceled. Flow balance constraints for each cargo is expressed as constraints (\ref{ConnCargoBalance1}) -- (\ref{ConnCargoBalance3}). It also restricts that for each cargo order, the un-shipped volume should be canceled. 
	The total volume shipped on each flight must not exceed the operating aircraft's capacity, according to constraint (\ref{ConnLegCap}).
	The set of short through connection constraints (\ref{ConnSC}) guarantees that each short through connection chosen in a cargo itinerary will be operated by a specific aircraft.
	
	One shortcoming of the arc-based model is the hardness to decide the copy size of flights, which lead to a vast network and solution inefficiency of the recovery problem.
	 To obtain a good solution more quickly, existing literature has proposed some pricing strategy to generate flight copies in a more smart way \citep{Liang2018,Huang2021}, which will improve the solution effectiveness significantly. Therefore, we also present a string-based model on the flight connection network and propose a column-and-row approach to solve the problem.
	
	\subsection{String-based Model (ACRP-S)}
	We rename the string-based model based on the flight connection network as ACRP-S. We introduce $L$ as the set of aircraft strings and $R$ as the set of cargo itineraries. We further introduce the set of departure time options for flight $f$ as $T_f$, and flight $f$'s duplication with departure time $t$ as $f_t$. Each short through connection in the string-based model is denoted as $(f_{t_1},f_{t_2}')$, with previous flight $f$ departs at $t_1$ and next flight $f'$ departs at $t_2$.
	The cost for each aircraft  $a$ that flies string $l$ is denoted as $c_{a,l}$ and the cost per unit of cargo $o$ shipped via itinerary $r$ is $c_{o,r}$. The cost of aircraft routing includes costs of aircraft swaps and delays for each flight composed in the string. The cost of cargo itinerary includes the cost of every flight change as well as the cost of the last flight's cargo delay if the cargo is delayed past its scheduled delivery time.
	The decision variables include whether to choose string $l$ for aircraft $a$ (denoted as $x_{a,l}$), and the volume of cargo $o$ shipped on itinerary $r$ (denoted as $w_{o,r}$).
	The additional notations for model ACRP-S are listed in Table \ref{Table:StringNotation}.
	
	\begin{table}[htbp] 
		\centering
		\small
		\caption{\label{Table:StringNotation}Additional notations for model ACRP-S} 
		\begin{tabular}{p{5em}p{30em}} 
			\toprule 
			\textbf{Sets} &\\
			$L   $& Set of aircraft strings indexed by $l$;\\
			$L_f $& Subset of aircraft strings that includes flight $f$;\\
			$R   $& Set of cargo itineraries indexed by $r$;\\
			$ T_f      $& Departure time set of flight $f$, indexed by $t$;  flight $f$'s original scheduled departure time is denoted as $t_f^0$; 
					flight $f$'s delay duration with departure time $t$ is denoted as $t_f^d$ equals to $t - t_f^0$; 
					duplication of flight $f$ with departure time $t$ is denoted as $f_t$;\\
			$SC  $& Set of short through connections indexed by $(f_{t_1},f_{t_2}')$, previous flight $f$ with departure time $t_1$ and next flight $f'$ with departure time $t_2$;\\
			
			\midrule
			\textbf{Parameters} &\\
			$c_{a,l} $& Cost incurred if aircraft $a$ flies string $l$;\\
			$c_{o,r}    $& Cost incurred for every unit of cargo $o$ with itinerary $r$;\\
			
			
			\midrule
			\textbf{Variables}&\\	
			$x_{a,l}$	& $x_{a,l} = 1$, if aircraft $a$ flies string $l$; otherwise, $x_{a,l}=0$;\\
			$w_{o,r}$	& Amount of cargo $o$ on itinerary $r$;\\
			
			\bottomrule 
		\end{tabular} 
	\end{table}
	
	Given the above notations, we present the  mathematical formulation of model ACRP-S as follows:

	\begin{align}
		\min &\sum_{f \in F} c_{f} y_{f}+\sum_{a \in A} \sum_{l \in L} c_{a, l} x_{a, l}+\sum_{o \in O} c_{o} z_{o}+ \sum_{o \in O}\sum_{r \in R} c_{o,r} w_{o,r}
		\label{MpObj}
		\\
		s.t. &\sum_{a \in A} \sum_{l\owns f} x_{a, l}+y_{f}=1, 
		\ \forall f \in F
		\label{MpCover}
		\\
		&\sum_{l \in L} x_{a, l} \leq 1,
		\ \forall a \in A
		\label{MpAft}
		\\
		& \sum_{r\in R}w_{o,r} + z_o = d_o, \forall o\in O \label{MpCargoCover}\\ 
		&\sum_{a \in A} \sum_{{l \owns f_t}} Cap_{a} x_{a, l} - \sum_{o\in O}\sum_{r\owns f_t}w_{o,r}\geq 0,
		\ \forall f \in F,\forall t \in T_f
		\label{MpCap} 
		\\
		&\sum_{a \in A} \sum_{{l \owns (f_{t_1},f_{t_2}')}} Cap_{a} x_{a, l} - \sum_{o\in O}\sum_{r\owns (f_{t_1},f_{t_2}')}w_{o,r}\geq 0,
		\ \forall (f_{t_1},f_{t_2}') \in SC  \label{MpShortConnect}\\ 
		&x_{a, l} \in\{0,1\}, 			\ \forall a \in A, \forall l \in L
		\label{MpX}\\
		&y_f \in\{0,1\},				\ \forall f \in F
		\label{MpY}\\
		&z_{o}\in \mathbb{Z}, 0 \leq z_{o} \leq d_o, 	\ \forall o \in O
		\label{MpZ}\\
		&w_{o,r} \in \mathbb{Z}, \forall o \in O, \forall r \in R
	\end{align}
	
	The objective function (\ref{MpObj}), like function (\ref{ConnObj}), is to minimize overall recovery costs, which are calculated as the summation of flight cancellation costs, aircraft routing costs, cargo cancellation costs, and cargo itinerary costs.  Similar to constraint (\ref{ConnFltCover}), constraint (\ref{MpCover}) is a set of constraints on cover for every flight.
	Every aircraft $a$ is restricted to flying no more than one string by constraint (\ref{MpAft}).
	For each cargo $o$, the volume that is not shipped by an itinerary must be canceled, according to constraint (\ref{MpCargoCover}).
	Constraints (\ref{MpCap}) -- (\ref{MpShortConnect}) are similar to (\ref{ConnLegCap}) -- (\ref{ConnSC}) that restrict the availability of aircraft capacity for each flight and each short through connection respectively.
	
	The model $\text{ACRP-S}$ is hard to be solved directly because the number of total possible aircraft strings and cargo itineraries is enormous.
	We present the solution approach in the next section.

	\section{Solution Approach for $\text{ACRP-S}$} \label{Sec:ColRow}
	In this section, we first propose a set of valid inequality constraints to improve the linear relaxation feasible space of the model. Then we describe the column-and-row generation framework and the details of aircraft string and cargo itinerary generation algorithms. For the purpose of accelerating the solution process, we further introduce a machine learning integrated algorithm to promote the selection of promising flight delay decisions.
	
	\subsection{Valid Inequality Constraints for ACRP-S} \label{Sec:ValidIneq}
\begin{figure}[htbp] 
		\centering 
		\includegraphics[width=0.6\textwidth]{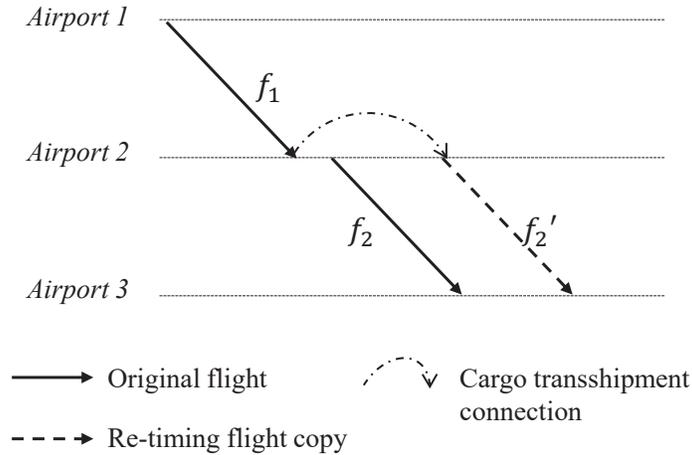} 
		\caption{An illustrative example} 
		\label{Figure_Example} 
	\end{figure}

	When the aircraft capacity is much larger than the cargo volume, the above model will provide a weak linear relaxation solution. To illustrate the situation, we provide the following example in Figure \ref{Figure_Example}.
	This example shows a recovery plan with one aircraft that operates two flights, which is $f_1$ and $f_2$, and one cargo needs to be shipped with the two flights. To ship the cargo, flight $f_2$ needs to be delayed to $f_2'$. The details of aircraft capacity, cargo demand, as well as feasible aircraft strings and cargo itineraries are shown in Table \ref{Table:VI_Example}.
	
	\begin{table}[htbp]
		\centering
		\caption{\label{Table:VI_Example} Details of the illustrative example}
		\begin{tabular} 
			{ll}
			\toprule
			Aircraft Capacity:
			& 15
			\\
			\midrule
			Aircraft Strings:
			& $x_1$: $f_1\rightarrow f_2$, cost equals to 0\\
			& $x_2$: $f_1\rightarrow f_2'$, cost equals to 50
			\\
			\midrule
			Cargo Demand:
			& 5
			\\
			\midrule
			Cargo Itineraries:
			& $w_1$: $f_1\rightarrow f_2'$, cost equals to 0
			\\
			\bottomrule
		\end{tabular}
	\end{table}%

	For convenience, we suppose all flights and cargoes are covered.
	Model $\text{ACRP-S}$ in this example can be written as:
	$\{	\min\ 50 x_2| x_1 + x_2 = 1; 15 x_2 \ge w_1; w_1 = 5\}$. 
	Here the irrelevant constraints are omitted. It is obvious that the optimal MIP solution for the aircraft string is $x_2^* =1$ with all volume of cargo covered, and the optimal objective value is $50$. However, the optimal solution for relaxed-LP is $x_1^*=2/3,x_2^*=1/3$ with an objective of $50/3$. It shows a great gap between the relaxed-LP and MIP results.
	Intuitively, the fractional solution only delays a fractional of an aircraft just enough to satisfy the demand.
	However, in the integer solution, if any fraction of a flight is delayed, the entire flight should be delayed.

	Therefore, to tighten the bounding constraints (\ref{MpCap}), we propose a set of additional valid inequalities illustrated as follows:
	
	\begin{align}
		\sum_{a \in A} \sum_{{l \in L_f \setminus L_f^t}} x_{a, l} + \frac{\sum_{r\owns f_t}w_{o,r}}{d_o}\leq 1, 
		\ \forall o \in O, \forall f \in F,\forall t \in T_f  
		\label{MpValidIneq} 
	\end{align}
	
	Here $L_f^t $ is defined as a subset of aircraft strings that includes flight $f$ with departure time $t$. In the constraint (\ref{MpValidIneq}), $(\sum_{r\owns f_t}w_{o,r})/{d_o}$ is an approximation for the fraction of aircraft strings that cover $f_t$ with cargo itineraries, and $\sum_{a \in A} \sum_{{l \in L_f \setminus L_f^t}} x_{a, l}$ computes the fraction of flight that does not depart on time $t$. Thus, the constraint limits the selection of multiple departure times for aircraft.
	It provides a good approximation when the total volume of the cargo (i.e., $d_o$) is small. After adding the valid inequality constraint, the feasible set of the LP relaxation of the MIP is improved.
	Continued with the example illustrated above, the valid inequality constraint in this case is $x_1 + w_1/5 \le 1$. The relaxed-LP solution is $x_2^* =1$, the same as the integer solution.
	For convenience, we refer to model ACRP-S with valid inequality constraints as $\text{ACRP-S}^*$ in the text that follows.
	
	\subsection{Column-and-Row Generation Framework} \label{subSec:ColRow}
	
	Column generation is a commonly used strategy for solving string-based models.
	By utilizing it, one can handle the enormous number of variables efficiently without generating all possible strings.
	However, column-generation can not be directly used to tackle the $\text{ACRP-S}^*$. 
	Specifically, when we generate new aircraft strings and cargo itineraries, we could obtain new flight delay decisions. 
	These new flight delay decisions should satisfy the capacity, short through connection and valid inequality constraints in (\ref{MpCap}), (\ref{MpShortConnect}) and (\ref{MpValidIneq}).
	 Because the new generated flight re-timing decisions can not be obtained in advance, it is impossible to ensure the solution is feasible when adding new delayed columns without adding additional associated constraints.
	
\begin{figure}[htbp] 
		\centering 
		\includegraphics[width=0.9\textwidth]{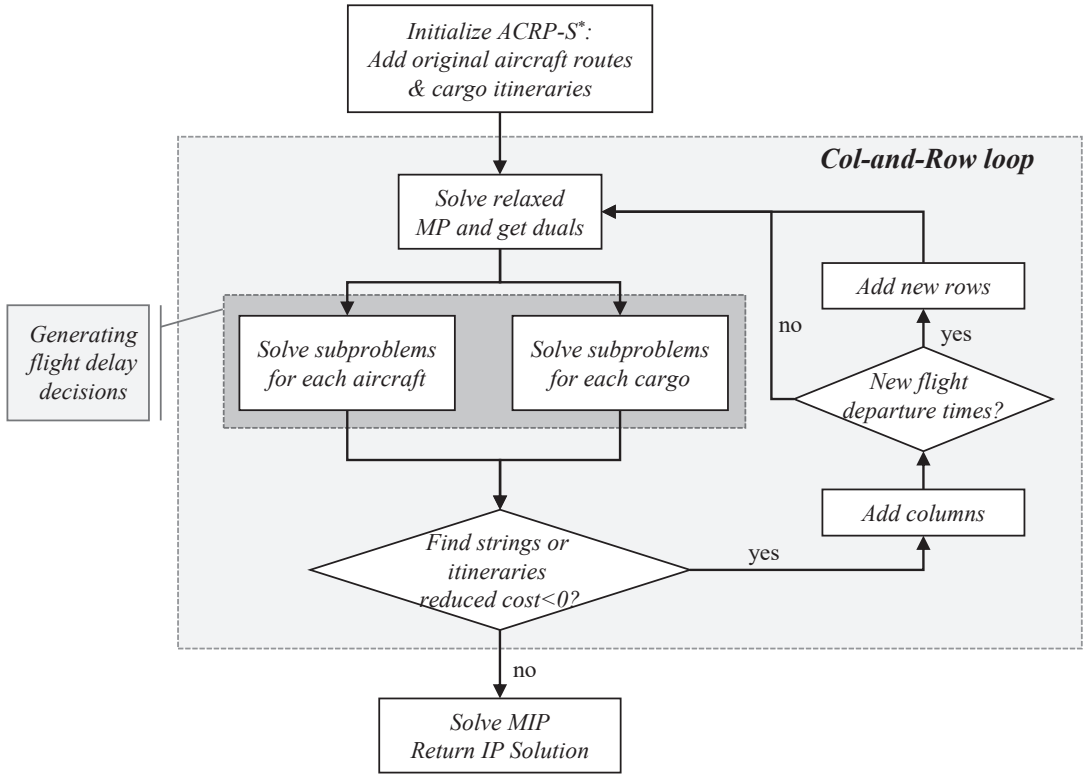} 
		\caption{Flow chart of the column-and-row generation framework for solving $\text{ACRP-S}^*$ } 
		\label{FigureFlowChart} 
	\end{figure}

	Thus we introduce a column-and-row generation approach to solve $\text{ACRP-S}^*$.
	This approach decomposes the original problem into three parts: 1) the master problem (i.e., the relaxed linear form of $\text{ACRP-S}^*$, denoted as MP) for selecting aircraft strings and cargo itineraries that minimize total recovery cost; 2) the sub-problems for generating better strings and itineraries; 3) the row generation process to ensure the feasibility of the new strings and itineraries. 
	The problems are solved iteratively until the linear relaxed master problem is optimal.
	Finally, we solve the mixed integer MP problem to obtain the integer solution. 
	The flowchart of the column-and-row generation approach in this study is shown in Figure \ref{FigureFlowChart}. 
		
	In detail, we solve the relaxed master problem to obtain dual values associated with each constraint for each iteration. The dual values are used as inputs for the sub-problems and serve as node (flight) weights to guide the new column generation process for each aircraft and cargo. 
	Better strings and itineraries with negative reduced costs generated in sub-problems are fed into the master problems until the optimal solution for the relaxed linear MP is obtained.
	Each time we add a new aircraft string or a new cargo itinerary (column) to the mater problem, we check the departure time for each flight included in the new aircraft string or cargo itinerary. If a new departure time and arrival time for a flight is generated, the master problem is updated with a set of constraints (\ref{MpCap}) and (\ref{MpValidIneq}) for the new departure time. Similarly, if a new short through connection is generated, a new constraint (\ref{MpShortConnect}) is added into the MP. 

	For convenience purpose, we denote the column-and-row generation solution approach described above as \textit{CRG}.
	In the following part of this section, we present details of the sub-problem solving algorithms used to generate better aircraft strings and cargo itineraries.
	
	\subsection{Sub-problem for Aircraft String Generation}
	
	\subsubsection{Reduced Cost Calculation for Aircraft String} \label{SecAftRC}

	In the aircraft string generation sub-problem, we first calculate the reduced cost for each string. For the relaxed master problem, we assume that $\alpha_f, \beta_a, \gamma_{f,t}, \eta_{(f_{t_1},f_{t_2}')}$ and $\pi_{o,f,t}$ are the dual variables associated to the constraints (\ref{MpCover}), (\ref{MpAft}), (\ref{MpCap}), (\ref{MpShortConnect}) and (\ref{MpValidIneq}) respectively. 
	
	Given aircraft $a$, the reduced cost $\bar{c}_{a,l}$ of its string $l$ is defined by Eq. (\ref{AftDual}).
	\begin{align}
		\bar{c}_{a,l}&=c_{a,l}-\sum_{f\in l}\alpha_f -\beta_a - \sum_{f_t\in l}Cap_a \gamma_{f,t} - \sum_{(f_{t_1},f_{t_2}')\in l}Cap_a \eta_{(f_{t_1},f_{t_2}')}\nonumber \\
		& - \sum_{f_t\in l}\sum_{o\in O}\sum_{t'\in T_f\setminus t}\pi_{o,f,t'}
		\label{AftDual}
	\end{align}
	
	The string cost $c_{a,l}$ of aircraft $a$ flying string $l$ is the summation of all the costs associated with swapping and delay cost of the flights in the string, i.e., $c_{a,l} = \sum_{f\in l}(c_{a,f}^{swap}+c_{f_t}^{delay})$. For each flight $f$ to aircraft $a$, the swap cost $c_{a, f}^{swap}>0$ if flight $f$ is not originally assigned to aircraft $a$, and $c_{a, f}^{swap}=0$ otherwise. 
	If flight $f$ is delayed, delay cost $c_{f_t}^{delay}>0$ and is linearly related to the delay duration in this study. 
	Thus, Eq.(\ref{AftDual}) can be rewrite as followed:
	
	\begin{align}
		\bar{c}_{a,l}&= - \beta_a + \sum_{f\in l}(c_{a,f}^{swap} -\alpha_f) - \sum_{f_t\in l}(Cap_a \gamma_{f,t} + \sum_{o\in O}\sum_{t'\in T_f\setminus t}\pi_{o,f,t'}- c_{f_t}^{delay})  \nonumber\\
		& \ \ \ \ -\sum_{(f_{t_1},f_{t_2}')\in l}Cap_a \eta_{(f_{t_1},f_{t_2}')} \label{AftDualReform}
	\end{align}
	
	The objective of the sub-problem is to get aircraft strings with negative reduced cost, or $\bar{c}_{a, l}<0$ equivalently.
	Eq.(\ref{AftDualReform}) shows that the reduced cost is composed of several parts. For a specific aircraft $a$, $\beta_a$ is independent of the flights assigned to the aircraft, $c_{a,f}^{swap} -\alpha_f$ are costs associated with each flight $f$ assigned to the aircraft. While the last part is related to the flights with respective departure time. In the column-and-row generation process, the last part is also determined by whether the constraints related to the departure time has been included in the master problem. If the constraints are not included in one iteration, this portion of the reduced cost is zero and will be updated in the following \textit{CRG} iterations.
	
	\subsubsection{Sub-problem Solution Algorithm for Aircraft String Generation} \label{SecAftSP}
	
	In the sub-problem for each aircraft, the goal is to find better strings with negative reduced costs that can potentially improve the optimal solution of the master problem. Each string is a sequence of flights from the given original airport to available destination airports for the aircraft. 

\begin{figure}[htbp] 
		\centering 
		\includegraphics[width=0.5\textwidth]{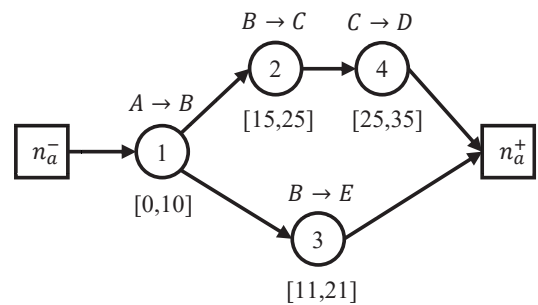} 
		\caption{An example of a flight connection network} 
		\label{FigureSPNetwork} 
	\end{figure}

	To solve the problem, we build a flight connection network $G_a(V_a,E_a)$ for each aircraft $a$ based on the original flight schedule. In the network, nodes represent flights with the originally scheduled origin, destination as well as the scheduled departure and arrival times. Each pair of flights $(i,j)$ is connected by an arc if flight $i$’s destination is the same as flight $j$’s origin, and flight $j$ departs later than flight $i$ in the original schedule. 
	We do not require flight $j$ to depart later than flight $i$'s arrival time plus turn time to meet connection availability, because the delay time is also a determination variable in the sub-problem. That is, the follow-up flight $j$ can be delayed until connected to flight $i$ with the connection requirement satisfied. This is different from the connection arc generation in model ACRP-A, in which the connection requirement must be satisfied when building connection arcs.
	A dummy source node $n_a^-$ and sink node $n_a^+$ are also added to the network. The source node represents the airport where and when the aircraft becomes available, and the sink node represents the end of the recovery horizon. 
	An illustration of a flight connection network for an aircraft departing from Airport $A$ is shown in Figure {\ref{FigureSPNetwork}}. Nodes in the shape of circles and squares represent flights and dummy nodes respectively. The origination and destination of each flight are shown above the nodes. Numbers in the brackets blow each node indicate the scheduled departure and arrival times (converted to integers) of each flight.

	The key decisions in the sub-problem include the set of flights to be assigned to the aircraft, the delay of each flight, as well as the connections between the flights. The goal is to find a path with negative reduced cost from the source node to the sink node. 
	We suggest using a multi-label shortest-path algorithm to solve the sub-problem. In addition to the reduced cost computed in \cref{SecAftRC}, we include flight delay time as an additional label because the delay for each flight is a recovery decision that affects connection feasibility for successive flights. 
	We use the shortest-path searching technique to decide on delays in accordance with earlier work by \cite{Liang2018}.
	For each flight, the delay decision depends on both the delay itself and the delay of its predecessors, and should be limited to the allowed delay range. More specifically, for two consecutive flights $(i, j)$, the delay of node $j$ should consider the arrival time with delay for node $i$ as well as the minimum turn time on the airport to satisfy the connectivity requirement. Therefore, the delay of each node contains propagated delay of its predecessor. 
	 Formally, we denote the label set of each node $i$ as $B_{a,i}$, in which each element is denoted as $b_{a,i}: \langle\bar{c}_{a,i}, t_{i}^{d}\rangle$.
	 In each label, $\bar{c}_{a,i}$ denotes the total reduced cost from the source node to node $i$ and $t_{i}^{d}$ denotes the delay duration of node $i$ itself.
	
	\begin{algorithm}[htb] 
		\caption{ Multi-label shortest path algorithm for aircraft string generation } 
		\label{Alg:AftSP} 
		\begin{algorithmic}[1] 
			\REQUIRE ~~\\ 
			Sort all the nodes (flights) in the chronological order.\\
			A flight connection network for aircraft $a: G_a(V_a, E_a)$, duals: $\beta_{a}$, $\alpha_f$, $\gamma_{f, t}$, $\pi_{o,f,t}$, $\eta_{(f_{t_1},f_{t_2}')}$.
			\ENSURE ~~\\ 
			String set $L_a$ for aircraft $a$ with negative reduced cost. 
			\STATE Set label set of the source node $B_{a,n_a^-}$ as $\{\langle 0,0\rangle\}$ and label sets of each other node as $\emptyset$.
			\FOR {node $i\in V_a$}
			\FOR {node $j \in Adj[i]$}
			\STATE $\quad$ Process arc $(i, j)$ (see \textbf{Algorithm \ref{Alg:Process Flight Arc}})
			\ENDFOR
			\ENDFOR
			\STATE Select non dominated label set $B_{a,n_a^+}^*$ 
			\FOR {lable $b_{a,n_a^+} \in B_{a,n_a^+}^*$ }
			\STATE $\bar{c}_{a,n_a^+} = \bar{c}_{a,n_a^+}-\beta_a$
			\IF{$\bar{c}_{a,n_a^+} < 0$ }
			\STATE Construct a new aircraft string $l_a$ by tracing back the predecessors of $b_{a,n_a^+}$
			\STATE $L_a = L_a \cup \{l_a\}$
			\STATE Add new generated departure times to each flight's departure time set 
			\ENDIF
			\ENDFOR
			\RETURN Aircraft string set $L_a$
		\end{algorithmic}
	\end{algorithm}
	
	Algorithm \ref{Alg:AftSP} presents a pseudo-code summarizing the multi-label shortest path algorithms, and the details of the arc processing algorithm are detailed in Appendix (see Algorithm \ref{Alg:Process Flight Arc}).
	As demonstrated in Algorithm \ref{Alg:AftSP}, the shortest path algorithm is initialized with the label set of the source node set as $\{\langle 0, 0\rangle\}$ and label sets of other nodes as $\emptyset$. 
	Then nodes in the network are sorted chronologically by departure time. 
	Once a processor node $i$ has been checked, we check each successor node $j$ by processing arc $(i, j)$. 
	The arc processing algorithm decides the delay duration for node $j$  (i.e. $t_j^d$) by the minimal delay time necessary to satisfy connection availability, or the arrival time plus the minimum turn time for the sequential flights. Formally, for each $t_j$ in flight $j$'s rescheduled departure time set, the delay duration of flight $j$ is $t_j^d=max\{t_i + fly\_time_i + turn\_time_{ij} , t_j \} - t_j^0$. The total path cost $\bar{c}_{a,j}$, from the source node to node $j$, is also obtained as $\bar{c}_{a,j} = \bar{c}_{a,i} + c_{aij}$, in which $c_{aij}$ consists of swap cost if flight $j$ is not originally assigned to aircraft $a$, delay cost of flight $j$ based on the delay duration, and dual for flight $j$ with rescheduled departure time. 
	
	Once we get a new label for node $j$, we check whether the delay time is feasible according to the maximum allowed delay time for each flight, as well as the aircraft's available time window.
	For node $i$, label $l_1$ dominates label $l_2$ if and only if two conditions are satisfied: 1) $\bar{c}_{a,j,1} \le \bar{c}_{a,j,2}$; and 2) $t_{j,1}^d \le t_{j,2}^d$.
	If the new label is not dominated by any other existing label, it is added to the node $j$'s label set.
	After processing every arc leading to the sink node, we compute the reduced cost for each label using the formula $\bar{c}_{a, l}=\bar{c}_{a,n_a^+}-\beta_{a}$, and obtain the set of superior labels with negative reduced cost. For each label in the superior label set, we trace back the label to get the predecessors and return better new generated aircraft strings to the master problem.

	\subsection{Sub-problem for Cargo Itinerary Generation} \label{Sec:CargoSP}
	
	\subsubsection{Reduced Cost Calculation for Cargo Itinerary} \label{SecCargoSP}
	The sub-problem for cargo itinerary generation is similar to the aircraft string generation process.
	We assume $\theta_o$ as the dual variable for constraint (\ref{MpCargoCover}). Then the reduced cost for each cargo itinerary is: 
	\begin{align}
		\bar{c}_{o,r}&=c_{o,r}+\sum_{f_t\in r} \gamma_{f,t}-\theta_o + \sum_{(f_{t_1},f_{t_2}')\in r} \eta_{(f_{t_1},f_{t_2}')} - \sum_{f_t\in r}\pi_{o,f,t}/d_o 
		\label{CargoDual}
	\end{align}
	
	The cargo itinerary cost $c_{o,r}$ of each unit of cargo $o$ shipped on itinerary $r$ is the summation of all the costs associated with flight change of each flight and the delay of the last flight.
	The delay cost for cargo occurs on the last flight in the cargo itinerary, indicated as $c_{o,r}^{delay}$.
	Therefore, $c_{o,r} = \sum_{f\in r}c_{o,f}^{change} + c_{o,r}^{delay}$. Eq.(\ref{CargoDual}) can thus be rewritten as follows:
	
	\begin{align}
		\bar{c}_{o,r}&= -\theta_o + \sum_{f\in r}c_{o,f}^{change} + c_{o,r}^{delay} + \sum_{f_t\in r} (\gamma_{f,t} - \pi_{o,f,t}/d_o) + \sum_{(f_{t_1},f_{t_2}')\in r} \eta_{(f_{t_1},f_{t_2}')}
		\label{CargoDualReform}
	\end{align}
	
	Eq.(\ref{CargoDualReform}) demonstrates that the reduced cost is made up of several components. 
	For a given cargo $o$, $\theta_o$ is independent of the flights included in the itinerary, while the other components are related to the flights contained. Specifically, $c_{o,f}^{change}$ is related to each flight assigned regardless of the departure time. $c_{o,r}^{delay}$ is linearly related to the delay duration of last flight. The rest portion are related to the flights with specific departure times. If the departure times are newly generated, the portion equals to zero because the constraints have not been included in the master problem and will be updated in the following iterations.
	
	\subsubsection{Sub-problem Solution Algorithm for Cargo Itinerary Generation}
	Finding better cargo itineraries with reduced cost is the main objective of the sub-problem for cargo itinerary generation.
	To solve the problem, we also build a flight connection network $G_o(V_o,E_o)$ for each cargo $o$. 
	Then we add source node $n_o^-$ and sink node $n_o^+$ to the network, similar to what we did in the aircraft's network. 
	Source node $n_o^-$ is connected to flights that depart from the origination of the cargo with the originally scheduled departure time no earlier than the first flight in the original cargo itinerary. 
	The distinction with the aircraft string generation sub-problem is that only the flights arrive at the cargo's destination are connected to the sink node $n_o^+$ in the cargo's network.
	
	In this study, we also make flight delay decisions in the sub-problem for cargo itinerary generation to get better integrated recovery decisions.
	In ACRP, the recovery decision of flight delay has "double-edged" consequences for the overall recovery solution.
	\begin{itemize}
		\item Shorter flight delays result in smaller flight delay and cargo delay costs. Moreover, for flight connections in cargo itineraries, shorter delay for the previous flight also results in fewer connection disruptions too.
		\item Longer flight delays might show some advantage in many situations. Specifically, for successive flights, it might be better for successive flights to delay longer in order to wait for the transshipment cargoes.
	\end{itemize}
	
\begin{figure}[htbp] 
		\centering 
		\includegraphics[width=0.7\textwidth]{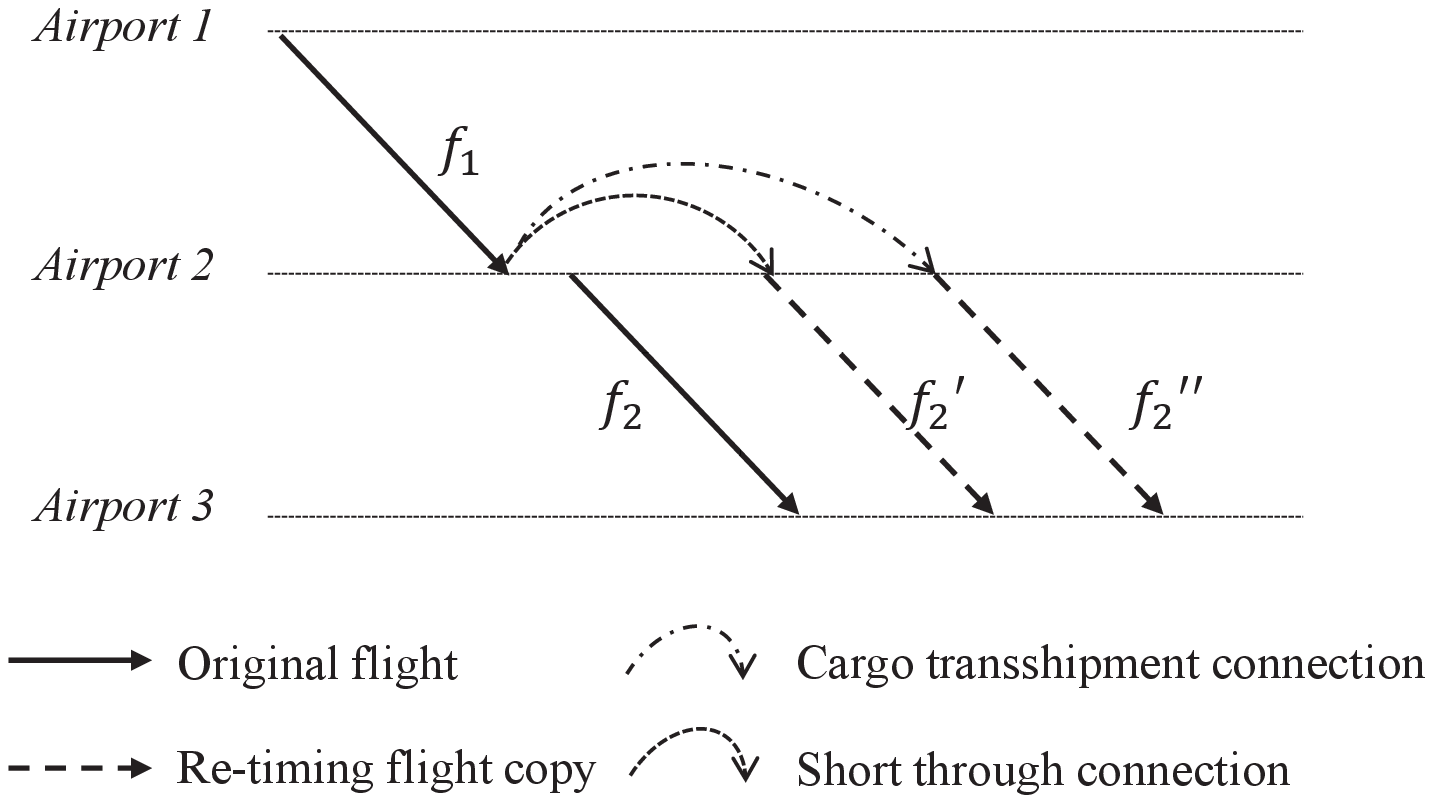} 
		\caption{An illustration of flight delay for cargo connections\\($f_2'$: re-timing flight with departure time that satisfies short through connection; $f_2''$: re-timing flight with departure time that satisfies standard transshipment.)} 
		\label{Fig:CargoDelay} 
	\end{figure}

	Therefore, different from the aircraft string generation algorithm, we generate two departure times for successive flights: one satisfies turn time of sequential flights operated by the same aircraft (short through connection), and the other satisfies the standard transshipment time which is sufficient for ground operations between flights operated by different aircraft. 
	As previously mentioned, taking advantage of short through connections can improve the efficiency of the recovery solution, and cargo can guide the generation of new short through connections.
	The flight delay decisions for cargo connections are illustrated in Figure \ref{Fig:CargoDelay}.
	
	We also apply the multi-label shortest-path algorithm to generate better cargo itineraries. The details are presented in Algorithm \ref{alg:CargoSP}.
	We denote the label set of each flight $i$ as $B_{o,i}$, and each element $b_{o,i}$ in the label set as $\left\langle\bar{c}_{o,i}, t_{i}^{d}, \{sc\}_i \right\rangle $, where $\bar{c}_{o,i}$ denotes the summation of cost from the source node to node $i$, $t_{i}^{d}$ denotes the delay duration on node $i$, and $\{sc\}_i$ denotes the set of short through connections from the source node to node $i$. 
	After sorting the flight nodes in chronological order, we process each flight connection arc $(i,j)$ to obtain labels for flight $j$ based on the existing labels for flight $i$. For each label of flight $j$, we generate a time set $ T_j^{new}$ consisting two departure times, which is $\{t_{i} + fly\_time_i+ turn\_time_{ij}, t_{i} + fly\_time_i+ trans\_time\}$. 
	We also have a departure time set for flight $j$ as $T_j$ that contains all the delay decisions generated in previous iterations. The details of the departure time set is shown in \cref{Sec:FlightRetimingSet}.
	We let $t_j^d = t_j-t_j^0$ for each $t_j \in T_j \cup T_j^{new}$ with $t_j \ge t_i + fly\_time_i +turn\_time_{ij}$.
	The total path cost $\bar{c}_{o,j}=\bar{c}_{o,i}+c_{oij}$, in which $c_{oij}$ consists of flight change of flight $j$ and duals related to flight $j$ with specific departure time. We set $\{sc\}_j = \{sc\}_i \cup \{(i,j)\}$ if the connection time for flight $i$ and $j$ with delay duration $t_i^d$ and $t_j^d$ is insufficient for standard transshipment, and $\{sc\}_j = \{sc\}_i$ otherwise.
	The pseudo-code of arc processing in the cargo network is detailed in Appendix A (see Algorithm \ref{Alg:Process Flight Arc_Cargo}). 
	After all arcs have been processed, we add $-\theta_o$ and the delay cost $c_{o,r}^{delay}$ to $\bar{c}_{o,n_o^+}$ and get the reduced cost for the itinerary. Then we select non-dominated non-negative labels from the candidate label set and add the corresponding better cargo itineraries to the master problem.
	
	It is noteworthy that we include  $\{sc\}_i$ in the label $b_{o,i}$ to prevent labels with short through connections from dominating other labels with standard transshipment time. The reason is that the delay duration for labels with short through connection is shorter, whereas the reduced cost $\bar{c}_{o,i}$ is usually overestimated when the dual values are missing due to delayed row generation. 
	More specifically, when we obtain a new departure time $t$ for flight $i$ during the column generation process, the fraction of the reduced cost related to the departure time $f_t$, i.e. $\sum_{f_t\in r} (\gamma_{f,t} - \pi_{o,f,t}/d_o) + \sum_{(f_{t_1},f_{t_2}')\in r} \eta_{(f_{t_1},f_{t_2}')}$, is missing.  
	For the dual problem, we have $\gamma_{f,t}\ge 0$, $\pi_{o,f,t}\le 0$ and $\eta_{(f_{t_1},f_{t_2}')}\ge 0$. Thus missing the new generated departure time overestimates the benefit of the new cargo itinerary. 
	Because the delay duration $t_i^d$ of the label consisting of standard transshipment is larger than labels consisting of short through connections, the former is thus dominated by the latter when both departure times are newly generated. 
	We add the component $\{sc\}_i$ into the label $b_{o,i}$ to deal with the issue. For node $i$,  label $l_1$ dominates label $l_2$ if and only if the following three conditions are satisfied: 1) $\bar{c}_{o,j,1} \le \bar{c}_{o,j,2}$; 2) $t_{j,1}^d \le t_{j,2}^d$; and 3) $\{sc\}_{i,1} \subseteq \{sc\}_{i,2}$.
	
	\begin{algorithm}[htb] 
		\caption{Multi-label shortest path algorithm for cargo itinerary generation} 
		\label{alg:CargoSP} 
		\begin{algorithmic}[1] 
			\REQUIRE ~~\\ 
			Sort all the nodes (flights) in the chronological order.\\
			A flight connection network for cargo $o: G_o(V_o,E_o)$ with the source node $n_o^-$ at the origination of cargo $o$ and sink node $n_o^+$ at the destination of cargo $o$, duals: $\theta_{o}$, $\gamma_{f, t}$, $\pi_{o,f,t}$, $\eta_{(f_{t_1},f_{t_2}')}$.
			\ENSURE ~~\\ 
			Itinerary set $R_o$ for cargo $o$ with negative reduced cost. 
			\STATE Set label set of the source node $B_{o,n_o^-}$ as $\{\langle 0,0,\{\} \rangle \}$ and label sets of each other node as $\emptyset$.
			\FOR {node $i\in V_o$}
			\IF {node $i$ arrives at the destination of cargo $o$}
			\STATE Continue loop
			\ENDIF
			\FOR {node $j\in Adj[i]$ }
			\STATE Process arc $(i,j)$ (see \textbf{Algorithm \ref{Alg:Process Flight Arc_Cargo}})
			\ENDFOR
			\ENDFOR
			\STATE Select non dominated label set $B_{o,n_o^+}^*$ 
			\FOR {lable $b_{o,n_o^+} \in B_{o,n_o^+}^*$}
			\STATE $\bar{c}_{o,n_o^+} = \bar{c}_{o,n_o^+} + c_{o,n_o^+}^{delay}-\theta_o$
			\IF{$\bar{c}_{o,n_o^+} < 0$ }
			\STATE Construct new cargo itinerary $r_o$ by tracing back the predecessors of $b_{o,n_o^+}$
			\STATE $R_o = R_o\cup \{r_o\}$
			\STATE Add new generated departure time to each flight's departure time set 
			\ENDIF
			\ENDFOR
			\RETURN Better cargo itinerary set $R_o$
		\end{algorithmic}
	\end{algorithm}

	\subsection{Flight Delay Decisions in Column-and-Row Generation} \label{Sec:FlightRetimingSet}
	As previously mentioned, we aim to make integrated flight delay recovery decisions that take into account both aircraft and cargo flight connection requirements.
	In each iteration of the column-and-row generation loop, however, the sub-problems for aircraft and cargoes are solved independently. 
	In order to include flight re-timing decisions generated by the aircraft string sub-problems into the cargo itinerary sub-problems or vice versa, we share the re-timing departure time set for each flight $f$ (denoted as $T_f$ in Table \ref{Table:StringNotation}) for all the aircraft and cargoes. 
	More specifically, when we process each arc $(i,j)$, we not only generate delay decisions for flight $j$ based on the connection requirement but also include all re-timing departure time decisions generated in previous iterations. At the end of each sub-problem, we check whether new flight departure times or new short through connections are generated. If this is the case, we update $T_f$ with the new generated flight departure times for each flight $f$.
	Thus, we created a "delayed" communication for the sub-problems with a shared re-scheduled flight departure time set.

	\subsection{Machine Learning Based Column-and-Row Generation} \label{Sec:MLTaskDefine}
	In the solution method described above, we incorporate as many potential flight delay decisions as possible into our model using a column-and-row generation strategy. 
	It is obvious that the model complexity and solution time are highly related to the number of flight delay decisions we introduced. That is, if we add too many flight delay decisions, the solution time for the master problem may become excessively long. We also noticed that a considerable fraction of the short through connection related flight delay decisions generated in the cargo itinerary sub-problems are not feasible for aircraft strings and thus are not selected in the LP solutions.
	
	Therefore, in order to reduce the number of "bad" flight delay decisions, we propose a machine learning approach to "smartly" select promising short through connections that are added to the master problem. The goal is to reduce the size of the problem and find effective solutions more quickly.
	In this study, we choose the popular decision tree model to select the promising short through connections and incorporate it into our column-and-row generation approach.
	
	\subsubsection{Prediction of Short Through Connections}
	The goal of the prediction is to identify whether a non-short through connection would then become a new short through connection during column-and-row generation.
	With the original flight schedule and all the cargo orders given, we can enumerate all the possible flight connections before the recovery process. 
	Specifically, we are interested in the connections that are not short through connections in the original schedule.
	We define $NSC$ as the set of non-short through connections in cargo itinerary networks.
	We formulate the prediction task as follows: for each connection $con\in NSC$, predict the short through connection class label of connection $con$, denoted as $p_{con}$.
	 This prediction task is a binary classification problem in which each connection is labeled "positive" if it becomes a newly generated short through connection in any iteration ($p_{con}=1$). On the other hand, the rest flight connections in $NSC$ are labeled "negative" ($p_{con}=0$).
	
	\subsubsection{Decision Tree Model} \label{SubSec:MLDecisionTree}
	A decision tree is a supervised classification procedure that recursively partitions a data set into smaller subdivisions.
	Decision trees are non-parametric and do not require assumptions regarding the distribution of the input data. This method can also handle nonlinear relations between features and classes.
	We used the classification and regression tree (CART) model described by \cite{breiman1984classification}. 
	
	In the decision tree algorithm, noise may cause some irrelevant features to be included among the selected tests and in turn causes \emph{overfitting}. Several pruning models have been proposed to prevent the fitting of noise \citep{Breslow1997simplifying}. In this research, we choose a post-pruning approach to correct the tree for overfitting. The decision tree is initially grown to its maximum size (denoted as $T_0$) and then trimmed from the bottom-up. We applied the Minimal Cost Complexity Pruning (MCCP) algorithm developed by \cite{breiman1984classification} that successively prunes the subtrees yielding minimal cost complexity. Rather than considering every possible subtree, MCCP generates a sequence of trees indexed by a non-negative tuning parameter $\alpha$. When $\alpha = 0$, the subtree simply equals to $T_0$. When we increase $\alpha$ from zero, we can get a sequence of nested smaller subtrees as a function of $\alpha$, ended with the root tree. 
	To select the best value of $\alpha$ as well as the corresponding best tree, we follow the algorithm using cross-validation described in \cite{James2013stanford}.
	
	It is notable that the machine learning data set in this study is imbalanced. This means that the number of instances representing the class of interest is significantly lower than that of the other class. The \textit{Imbalance Ratio} (IR) is a typical approach to describe the level of imbalance in a two-class problem, and is defined as the number of negative class instances divided by the number of positive class instances \citep{orriols2009}.
	The data set in this study is imbalanced because the number of "negative" (non-short through connections) instances are apparently larger than "positive" (short through connections) ones. If we treat the positive and negative instances equally while estimating the model performance, we would concentrate on the majority class and ignore the minority class \citep{fernandez2018}.
	However, it is more crucial to identify the minority class, i.e. potential short through connections, than it is to avoid misclassifying a connection as a short through connection in our study. 
	We can solve the issue by incorporating asymmetric misclassification costs into the model evaluation process.
	Cost $c(i|j)$ denotes a cost of misclassifying a class $j$ instance as belonging to class $i$, while  $c(i|j)=0$ for $i = j$.
	Without loss of generality, we impose a unity condition, at least one $c(i|j)=1$, which is the minimum misclassification cost. The unity condition allows us to measure the number of high cost errors. 
	We employ the most popular heuristic approach, which estimates the cost directly using the IR index. 
	In this set-up, $c(0|1)=IR$ and $c(1|0)=1$. 
	The cost matrix is then incorporated into the cost function of the pruning process of MCCP as described by \cite{breiman1984classification}.
	In detail, if we assign each leaf node $t$ to class $i$, the node's misclassification cost is $\sum_j c(i|j)p(j|t)$, in which $p(j|t)$ is the proportion of class $j$ in node $t$. Then the node is assigned to a class that minimizes the misclassification cost. The misclassification cost of note $t$ can be then defined as $r(t) = min_i\sum_j c(i|j)p(j|t)$. The total misclassification cost for the tree $T$ is $R(T) = \sum_{t\in\widetilde{T}}r(t)p(t)$, in which $\widetilde{T}$ denotes the set of all leaf nodes and $p(t)$ denotes the proportion of node $t$ in the tree.

	To conduct the decision tree model training, we first obtain training data from solutions of earlier days. Disruptions of earlier days are solved by the column-and-row generation process and solution details are recorded for training. 
	Since the flight network structure and cargo order characteristics are not likely to change dramatically, the solutions are usually suitable for training.
	The full process of the decision tree model training in this study is described as follows. 
	First, the total data set $D$ is partitioned into $D.train$ (for training) and $D.test$ (for testing). We generate the original full tree with $D.train$, and apply MCCP to obtain a sequence of best subtrees $\{T_{\alpha}\}$ and a corresponding discrete parameter set of $\alpha$ (denoted as $\mathcal{A}$).
	After that, we get the best tree with inner-run K-fold cross-validation. In detail, we divide $D.train$ into K folds. For each fold $k \in \{1,\cdots,K\}$, the $k$th fold $D.val_k$ is for validation and the rest $D.tr_k$ is for training. In each iteration of the inner-run, we apply the MCCP criterion to get the optimal decision tree for each $\alpha \in \mathcal{A}$ with data set $D.tr_k$. Then we test the misclassification cost on $D.val_k$. As a result, we can select a value of $\alpha^*$ that minimizes the average misclassification cost for the validation sets, and obtain the corresponding tree $T_{\alpha^*}$ as the output for the prediction task.
	At last, we can evaluate misclassification cost as the generalization performance of the generated decision tree on the test set $D.test$. 
	
	\subsubsection{Solution Approach for Machine Learning Based Column-and-Row Generation} \label{SubSec:CGML}
	Based on the column-and-row generation solution approach illustrated in Figure \ref{FigureFlowChart} (denoted as \textit{CRG}), we incorporated the above mentioned machine learning prediction process into the column-and-row generation approach. The integral solution approach is represented in Figure \ref{FigureFlowChart}. 
	Specifically, we train the short through connection prediction model first and save the decision tree model. Then we generate the non-short through connection set NSC before the column-and-row generation process. In the column-and-row generation iterations, after we get the LP solution of the master problem, we summarize the features for each connection and utilize the pre-trained decision tree to predict the short through connection class label for each connection $con$, which is denoted by $\hat{p}_{con}$. 
	The arc processing algorithm in the cargo itinerary generation sub-problems is modified with the short through connection result. If connection $(i, j)$ is not predicted as a short through connection ($\hat{p}_{con}=0$), we only generate the delay decision for flight $j$ that meets the standard transshipment connection time and ignore the delay decision that meets the short through connection time. In cases that $\hat{p}_{con}=1$, the arc processing algorithm is unchanged.
	For convenience purpose, the solution approach of integral column-and-row generation with prediction is denoted as \textit{ML-CRG} in the following context.
	
	\begin{figure}[htbp] 
		\centering 
		\includegraphics[width=0.8\textwidth]{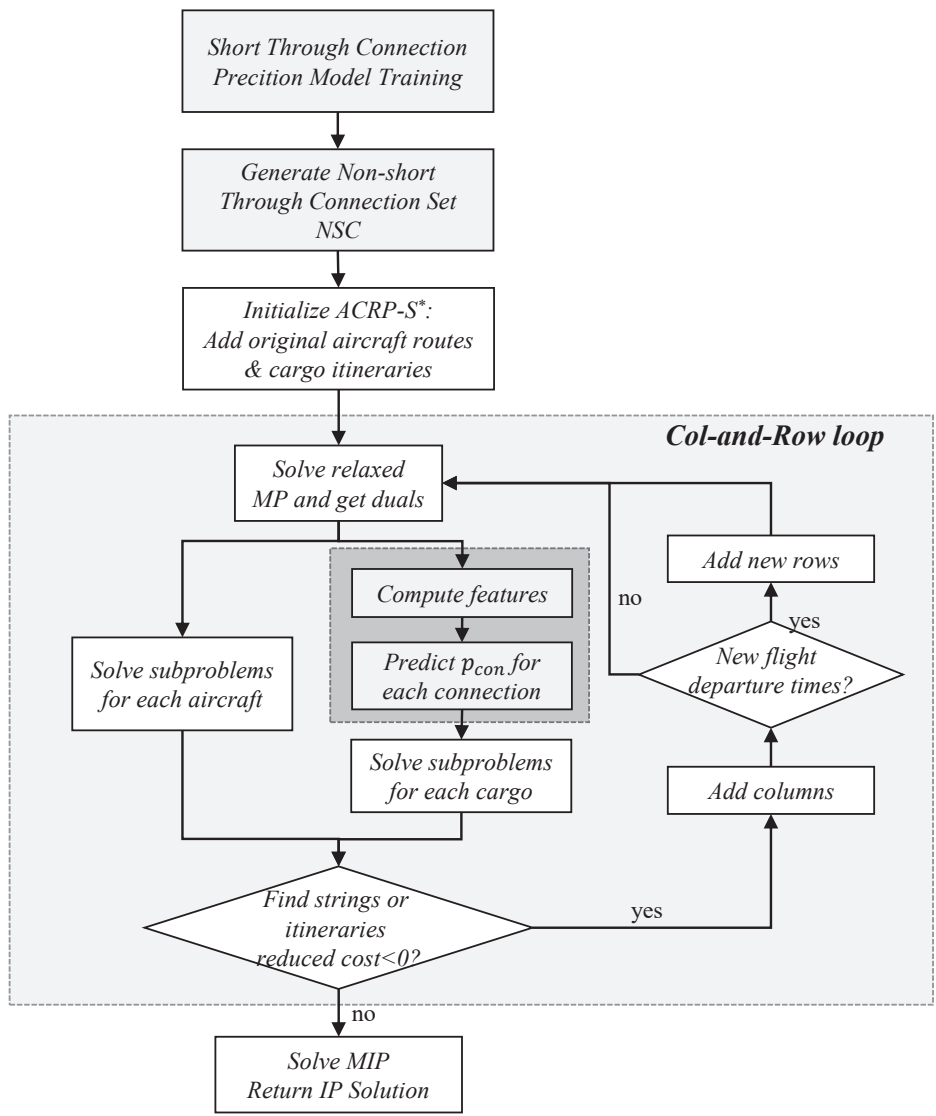} 
		\caption{Flow chart of the machine learning based column-and-row generation} 
		\label{FigureFlowChartML} 
	\end{figure}
	
	\section{Computational Study} \label{Sec:Experiments}
	
	In this section, to validate the effectiveness and solution quality of our proposed model and solution method, we conduct our computational testing for real-world recovery scenarios that affect freighter airlines.
	After describing our data and disruption scenarios, we compare the performance of the string-based integrated recovery model $\text{ACRP-S}^*$ with the arc-based model ACRP-A. 
	We also present the performance comparison between the integrated model with the sequential recovery model (denoted as ACRP-SEQ). The performance of the recovery solution quality is measured by use of recovery policies, overall recovery cost, and the computational performance is defined as the time consumed to solve the recovery problem.
	Finally, we present the full process of machine learning prediction as well as the effectiveness of the machine learning based column-and-row generation approach.
	
	\subsection{Description of Data and Disruption Scenarios} \label{Sec:ExperimentData}
	
	The data used in this study is based on two sets of real operational data from the largest Chinese freighter airlines that serves 46 airports. The freighter airline's AOC department makes recovery decisions within a time window of 1-2 days. The two data sets are labeled as F102-A34 and F284-A65, respectively. 
	The smaller data set F102-A34 is a 36-hour time window schedule operated by a single fleet Boeing 757-200. This schedule consists of 34 aircraft, 102 flights, and 242 cargo orders.
	The other schedule, F284-A65, is a two-day multi-fleet schedule with 65 aircraft, 284 flights and 551 cargo orders. Boeing 737-300, 757-200, 767-300, and 747-400 are among the fleets, with ULDs measuring their capacity ranging from 8 to 38. 
	In both data sets, each aircraft has a pre-scheduled string, and each cargo order has a pre-scheduled itinerary as well as a due time.
	
	We model a number of disruption scenarios using the original data, including a wide range of potential disruptions, such as flight time overlap (FO), airport mismatch (AM), capacity shortage (CS), aircraft on ground (AOG), closure of a hub airport (AC), and combinations of disruptions. 
	The characteristics of the scenarios are described in Table \ref{Table:Scenarios}. Scenarios 1-6 are generated from F102-A34, and Scenarios 7-12 are generated from F284-A65, respectively.

	\begin{table*}  
		\centering
		\caption{Scenarios Description} \label{Table:Scenarios}
		\begin{threeparttable}
			\begin{tabular}
				{lrrrrrl}
				\toprule
				Scenario & \# flt & \# aft & \# apt & cargo & Horizon (h) & Disruption Description \\
				\midrule
				1 & 102 & 34 & 33 & 242 & 36 & FO \& AM \& CS  \\
				2 & 102 & 34 & 33 & 242 & 36 & 48 h AOG  \\
				3 & 102 & 34 & 33 & 242 & 36 & 4 h AC  \\
				4 & 102 & 34 & 33 & 242 & 36 & FO \& AM \& CS \& 48 h AOG  \\
				5 & 102 & 34 & 33 & 242 & 36 & FO \& AM \& CS \& 4 h AC  \\
				6 & 102 & 34 & 33 & 242 & 36 & 48 h AOG \& 4 h AC  \\
				7 & 284 & 65 & 46 & 551 & 48 & FO \& AM \& CS  \\
				8 & 284 & 65 & 46 & 551 & 48 & 48 h AOG  \\
				9 & 284 & 65 & 46 & 551 & 48 & 4 h AC  \\
				10 & 284 & 65 & 46 & 551 & 48 & FO \& AM \& CS \& 48 h AOG  \\
				11 & 284 & 65 & 46 & 551 & 48 & FO \& AM \& CS \& 4 h AC  \\
				12 & 284 & 65 & 46 & 551 & 48 & 48 h AOG \& 4 h AC  \\
				\bottomrule
			\end{tabular}
		
			\begin{tablenotes}
				\footnotesize
				\item FO: flight time overlap; AM: airport mismatch; CS: capacity shortage; 
				\item AOG: aircraft on ground; 
				\item AC: airport closure
			\end{tablenotes}
		\end{threeparttable}
	\end{table*}

	Within the horizon, model ACRP-A and $\text{ACRP-S}^*$ implements a set of flight recovery policies, including flight cancellations, aircraft swaps, and fight delays, as well as cargo recovery policies, including cargo cancellations, flight changes, and cargo delays. Airlines incur costs for flight/cargo schedule changes in the recovery process. The costs are modeled in ACRP quantitatively as cost parameters for each recovery policy. The cost parameters for the recovery policies are given by airline companies. In practice, the cost parameters are determined by airline decision makers that show their preference for recovery policies. In all the scenarios, the disruption cost for flights and cargoes is shown in Table \ref{Table:Cost}. Note that each item of cost for a flight is charged for every flight, and the cost of cargo is penalized for every ULD unit of the cargo order.
	
	\begin{table}[htbp]
		\centering\small
		\caption{\label{Table:Cost}Cost Parameters of Test Scenarios}
		\begin{tabular}{ccccccc}
			\toprule
			\multicolumn{3}{c}{Flight related cost (per flight)} & & \multicolumn{3}{c}{Cargo related  cost (per ULD)}                                                                            \\
			\cmidrule{1-3} \cmidrule{5-7}
			\multicolumn{1}{c}{Cancellation} & 
			\multicolumn{1}{c}{\begin{tabular}[c]{@{}c@{}}Aircraft \\ swap\end{tabular}} & 
			\multicolumn{1}{c}{\begin{tabular}[c]{@{}c@{}}Delay cost \\ per hour\end{tabular}} & 
			& 
			\multicolumn{1}{c}{Cancellation}  &
			\multicolumn{1}{c}{\begin{tabular}[c]{@{}c@{}}Flight \\ change\end{tabular}} &
			\multicolumn{1}{c}{\begin{tabular}[c]{@{}c@{}}Delay cost \\ per hour\end{tabular}}
			\\
			\midrule
			1200  & 40 & 120  & & 60 & 1 & 2.4  \\
			\bottomrule                             
		\end{tabular}
	\end{table}

	For all scenarios, the short through connection time for cargo depends on the turn time for sequential flights operated by the same aircraft, varying from 70 to 90 minutes,  according to the fleet type of the aircraft, type of the flight (domestic or international) as well as the throughout of the airport. And the standard transshipment time for cargoes is 2 hours between different aircraft.
	
	All the experiments are implemented in C++ by calling CPLEX 12.7 as the linear programming and the mixed inter programming solver. The program was run on a server with a 2.90 Hz Intel Xeon 8268 CPU, 128G RAM, and a Windows Server 2019 system.
	
	\subsection{Computational results for model ACRP-A } \label{Sec:ExperimentArc}
	
	The scale of the arc-based model depends on the duplication size of every flight, which is determined by the maximum delay duration and the delay interval. However, it is hard to get the optimal delay duration for the disrupted schedule before making recovery decisions. The maximum delay duration is around 3 to 4 hours according to practical experience.
	The choice of delay interval is tactical because a smaller delay interval leads to a better solution with a longer solution time, whereas a larger delay interval leads to a shorter solution time and a worse solution quality.
	
	In this study, the delay interval for flight duplication is set as 5 minutes for small scenarios simulated from schedule F102-A34. The interval for scenarios of schedule F284-A65 is 30 minutes. 
	More precise delay decisions result in a problem size that is too large for our computational server to solve due to memory overflow.
	Furthermore, the maximum delay duration for scenarios of schedule F284-A65 has also been decreased to 3 hours due to computing capacity limitations.
	
	For the tested disruption scenarios, AOG is modeled as an additional flight that must be covered with an unchangeable timetable. And airport closure limits the feasibility of flight copies. The result of model ACRP-A for the above scenarios is presented in Table \ref{Table:ArcModelResult}. The total number of the variables, constraints and non-zero parameters are reported. Additionally, we report "Int. gap" as the gap between the best solution and the best bound found by the optimizing solver. 
	
	\begin{table}[htbp]
		\centering\small
		\caption{\label{Table:ArcModelResult}Computational Results of Model ACRP-A}
		\begin{tabular}{lrrrrrrrr}
			\toprule
			\multicolumn{1}{c}{Scenario}
			&\multicolumn{1}{c}{\begin{tabular}[c]{@{}c@{}}Max \\ delay\end{tabular}} &\multicolumn{1}{c}{\begin{tabular}[c]{@{}c@{}}Delay \\ interval\end{tabular}} 
			&\multicolumn{1}{c}{\# variables }
			&\multicolumn{1}{c}{\# constraints }
			&\multicolumn{1}{c}{\# non zeros }
			&\multicolumn{1}{c}{\begin{tabular}[c]{@{}c@{}}IP \\ Obj. \end{tabular}} 
			&\multicolumn{1}{c}{\begin{tabular}[c]{@{}c@{}}IP \\ time(s)\end{tabular}}
			& \multicolumn{1}{c}{\begin{tabular}[c]{@{}c@{}}Int. \\ gap\end{tabular} }
			\\
			\midrule
			1 & 240 & 5 & 71,043,931 & 797,473 & 244,622,732 & 632.20 & 803.0 & 0.00\%  \\
			2 & 240 & 5 & 70,572,563 & 795,900 & 243,155,829 & 738.40 & 3,610.0 & 0.00\%  \\
			3 & 240 & 5 & 70,572,061 & 795,899 & 243,157,416 & 1,604.60 & 753.0 & 0.00\%  \\
			4 & 240 & 5 & 71,044,433 & 797,475 & 244,624,736 & 1,229.80 & 3,076.0 & 0.00\%  \\
			5 & 240 & 5 & 71,043,931 & 797,474 & 244,626,323 & 2,110.80 & 790.0 & 0.00\%  \\
			6 & 240 & 5 & 70,572,563 & 795,901 & 243,159,420 & 2,210.00 & 2,642.0 & 0.00\%  \\
			7 & 180 & 30 & 48,747,976 & 788,261 & 163,934,304 & 1,754.40 & 4,831.0 & 0.00\%  \\
			8 & 180 & 30 & 48,805,924 & 788,640 & 164,068,416 & 2,499.20 & 8,906.0 & 0.00\%  \\
			9 & 180 & 30 & 48,805,914 & 788,639 & 164,243,992 & 2,805.20 & 3,320.0 & 0.00\%  \\
			10 & 180 & 30 & 48,747,986 & 788,264 & 164,109,818 & 5,175.60 & 9,891.0 & 0.00\%  \\
			11 & 180 & 30 & 48,747,976 & 788,262 & 164,109,784 & 3,615.60 & 5,005.0 & 0.00\%  \\
			12 & 180 & 30 & 48,805,924 & 788,641 & 164,244,026 & 4,365.20 & 7,909.0 & 0.00\%  \\
			\bottomrule                             
		\end{tabular}
	\end{table}
	
	We can see from Table \ref{Table:ArcModelResult} that half of the scenarios take longer than an hour to solve. The major reason for the low solving efficiency is that most of the duplicated flight delay decisions are not good for the recovery problem. 
	Contrarily, the algorithm we proposed for model $\text{ACRP-S}^*$ produces flight duplication when needed and can improve the solution efficiency, regardless of the extremely huge scale of potential aircraft strings and cargo itineraries.
	We present the computational results of model $\text{ACRP-S}^*$ next.

	\subsection{Computational results for model $\text{ACRP-S}^*$}
	In this section, we show the computational results for model $\text{ACRP-S}^*$ using algorithm \textit{CRG}.
	AOG is also modeled as an additional flight that must be covered. Flights that are interrupted by airport closure are delayed until the airport is available.
	
	We show the details of the column-and-row generation process and results in Table \ref{Table:StringResult}.	
	"LP Obj." is the final LP objective value of the column-and-row generation iteration until no better aircraft strings or cargo itineraries can be found. "IP Obj." is the IP solution obtained by solving the mixed-integer problem directly after the column-and-row iterations. 
	The "Int. gap" shows the gap between the LP optimal value and the final IP result, computed as (IP Obj. - LP Obj.)/LP Obj.. The small gap for all cases indicates that the relaxed LP result is a good approximation for the original mixed-integer problem.
	Additionally, we present the number of initial aircraft strings and the number of initial cargo itineraries as "Init \# x" and "Init \# w" separately. We also present details of the column-and-row generation process. "\# Col\ Iter" shows the iterations of the column generation, and "\# Row\ Iter" shows the iterations of the row generation. The iterations of row generation are fewer because we do not generate new flight departure times after column generation in every iteration. 
	The number of final aircraft strings and the number of cargo itineraries generated are shown as "Final \# x" and "Final \# w". "LP time(s)" is the run time (in seconds) for the column-and-row generating process, and "IP time(s)" is the total run time (in seconds) including the MIP solution.
	
	\begin{table}[htbp]
		\centering
		\caption{\label{Table:StringResult}Details of Column-and-Row Approach for Model $\text{ACRP-S}^*$}
		\resizebox{\linewidth}{!}{  
		\begin{tabular} 
			{lrrrrrrrrrrrr}
			\toprule
			Scenario 
			&\multicolumn{1}{c}{\begin{tabular}[c]{@{}c@{}} Init \\ \# x \end{tabular} } 
			&\multicolumn{1}{c}{\begin{tabular}[c]{@{}c@{}} Init \\ \# w \end{tabular} } 
			&\multicolumn{1}{c}{\begin{tabular}[c]{@{}c@{}} \# Col \\ Iter \end{tabular} } 
			&\multicolumn{1}{c}{\begin{tabular}[c]{@{}c@{}} \# Row \\ Iter \end{tabular} } 
			&\multicolumn{1}{c}{\begin{tabular}[c]{@{}c@{}} Final \\ \# x \end{tabular} } 
			&\multicolumn{1}{c}{\begin{tabular}[c]{@{}c@{}} Final \\ \# w \end{tabular} } 
			&\multicolumn{1}{c}{\begin{tabular}[c]{@{}c@{}} \# Add \\ Rows \end{tabular} } 
			
			&\multicolumn{1}{c}{\begin{tabular}[c]{@{}c@{}} LP \\ Obj. \end{tabular} } 
			&\multicolumn{1}{c}{\begin{tabular}[c]{@{}c@{}} LP \\ time(s) \end{tabular} } 
			&\multicolumn{1}{c}{\begin{tabular}[c]{@{}c@{}} IP \\ Obj. \end{tabular} } 
			&\multicolumn{1}{c}{\begin{tabular}[c]{@{}c@{}} IP \\ time(s) \end{tabular} } 
			&\multicolumn{1}{c}{\begin{tabular}[c]{@{}c@{}} Int. \\ gap \end{tabular} } 
			\\
			\midrule		
			1 & 34 & 242 & 6 & 4 & 457 & 273 & 457 & 594.86 & 2.8 & 632.20 & 3.4 & 6.28\%  \\
			2 & 34 & 242 & 9 & 8 & 761 & 339 & 583 & 706.33 & 2.7 & 738.40 & 2.9 & 4.54\%  \\
			3 & 34 & 242 & 8 & 6 & 952 & 363 & 680 & 1,572.59 & 3.2 & 1,604.60 & 3.5 & 2.04\%  \\
			4 & 34 & 242 & 9 & 8 & 907 & 359 & 647 & 1,192.39 & 3.7 & 1,229.80 & 3.9 & 3.14\%  \\
			5 & 34 & 242 & 7 & 4 & 987 & 365 & 676 & 2,073.45 & 2.8 & 2,110.80 & 3.0 & 1.80\%  \\
			6 & 34 & 242 & 8 & 6 & 1,149 & 419 & 798 & 2,177.92 & 2.8 & 2,210.00 & 3.1 & 1.47\%  \\
			7 & 65 & 551 & 37 & 35 & 9,822 & 3,381 & 6,368 & 1,310.21 & 719.5 & 1,359.80 & 727.4 & 3.78\%  \\
			8 & 65 & 551 & 14 & 13 & 13,076 & 2,232 & 5,136 & 2,221.94 & 311.0 & 2,256.60 & 314.8 & 1.56\%  \\
			9 & 65 & 551 & 25 & 21 & 10,266 & 4,240 & 9,128 & 2,455.42 & 1,415.5 & 2,495.00 & 1,424.0 & 1.61\%  \\
			10 & 65 & 551 & 33 & 30 & 10,471 & 4,107 & 7,990 & 4,658.83 & 1,119.4 & 4,713.40 & 1,131.1 & 1.17\%  \\
			11 & 65 & 551 & 32 & 28 & 9,458 & 3,805 & 7,105 & 3,098.87 & 908.9 & 3,153.40 & 918.6 & 1.76\%  \\
			12 & 65 & 551 & 27 & 19 & 9,662 & 2,968 & 6,624 & 4,015.39 & 660.5 & 4,055.00 & 665.5 & 0.99\%  \\
			\bottomrule
		\end{tabular}
	}
	\end{table}%

	In order to provide a more thorough view of the recovery solution, we present the recovery policy details for the test scenarios in Table \ref{Table:IPdetail}, and the breakdown of overall cost by each recovery policy in Table \ref{Table:IPcost}.
	
	\begin{table}[htbp]
		\centering\small
		\caption{\label{Table:IPdetail}Recovery Policy Details of Model $\text{ACRP-S}^*$'s IP Solutions}
		\begin{tabular}{lrrrrrrrr}
			\toprule
			\multicolumn{1}{c}{} & \multicolumn{3}{c}{Flight} & & \multicolumn{3}{c}{Cargo}                                                
			 \\
			\cmidrule{2-4} \cmidrule{6-8}
			\multirow{2}{*}{Scenario}  &                    
			\multicolumn{1}{c}{\begin{tabular}[c]{@{}c@{}}\# flight \\ cancel\end{tabular} } &
			\multicolumn{1}{c}{\begin{tabular}[c]{@{}c@{}}\# aircraft \\ swap \end{tabular}} & 
			\multicolumn{1}{c}{\begin{tabular}[c]{@{}c@{}} Delay \\ minutes\end{tabular} } & 
			& 
			\multicolumn{1}{c}{\begin{tabular}[c]{@{}c@{}}\# cargo \\ cancel\end{tabular} } &
			\multicolumn{1}{c}{\begin{tabular}[c]{@{}c@{}}\# flight \\ change\end{tabular} } &
			\multicolumn{1}{c}{\begin{tabular}[c]{@{}c@{}} Delay \\ minutes\end{tabular} } 
			\\
			\midrule
		
			1 & 0 & 3 & 75 &   & 5 & 4 & 1,455  \\
			2 & 0 & 3 & 240 &   & 1 & 1 & 1,935  \\
			3 & 0 & 0 & 700 &   & 1 & 1 & 3,590  \\
			4 & 0 & 6 & 280 &   & 5 & 4 & 3,145  \\
			5 & 0 & 3 & 740 &   & 5 & 4 & 5,170  \\
			6 & 0 & 3 & 905 &   & 1 & 1 & 5,475  \\
			7 & 0 & 5 & 425 &   & 2 & 27 & 4,070  \\
			8 & 1 & 5 & 160 &   & 8 & 33 & 590  \\
			9 & 0 & 5 & 970 &   & 2 & 33 & 5,050  \\
			10 & 1 & 5 & 1,235 &   & 8 & 27 & 8,410  \\
			11 & 0 & 5 & 1,235 &   & 2 & 27 & 8,410  \\
			12 & 1 & 5 & 970 &   & 8 & 33 & 5,050  \\
			\bottomrule                             
		\end{tabular}
	\end{table}

	\begin{table}[htbp]
		\centering\small
		\caption{\label{Table:IPcost}Recovery Cost Details of Model $\text{ACRP-S}^*$'s IP Solutions}
		\resizebox{\linewidth}{!}{  
		\begin{tabular}{lrrrrrrrrrr}
			\toprule
			\multicolumn{1}{c}{} & \multicolumn{4}{c}{Flight related costs} & & \multicolumn{4}{c}{Cargo related costs}  &                                                         \\
			\cmidrule{2-5} \cmidrule{7-10}
			\multirow{2}{*}{Scenario}                      
			& \multicolumn{1}{c}{\begin{tabular}[c]{@{}c@{}}Flight \\cancel cost\end{tabular} } 
			& \multicolumn{1}{c}{\begin{tabular}[c]{@{}c@{}}Aircraft \\swap cost\end{tabular}} 
			& \multicolumn{1}{c}{\begin{tabular}[c]{@{}c@{}}Flight \\delay cost\end{tabular} } 
			& \multicolumn{1}{c}{\begin{tabular}[c]{@{}c@{}}Total \\flight cost\end{tabular} }  
			&
			& \multicolumn{1}{c}{\begin{tabular}[c]{@{}c@{}}Cargo\\cancel cost\end{tabular} } 
			& \multicolumn{1}{c}{\begin{tabular}[c]{@{}c@{}}Flight\\change cost\end{tabular} } 
			& \multicolumn{1}{c}{\begin{tabular}[c]{@{}c@{}}Cargo\\delay cost\end{tabular} } 
			& \multicolumn{1}{c}{\begin{tabular}[c]{@{}c@{}}Total\\ cargo cost\end{tabular} }
			& \multicolumn{1}{c}{\begin{tabular}[c]{@{}c@{}}Recovery\\ cost\end{tabular} }
			\\
			\midrule
			1 & 0.00 & 120.00 & 150.00 & 270.00 &   & 300.00 & 4.00 & 58.20 & 362.20 & 632.20  \\
			2 & 0.00 & 120.00 & 480.00 & 600.00 &   & 60.00 & 1.00 & 77.40 & 138.40 & 738.40  \\
			3 & 0.00 & 0.00 & 1,400.00 & 1,400.00 &   & 60.00 & 1.00 & 143.60 & 204.60 & 1,604.60  \\
			4 & 0.00 & 240.00 & 560.00 & 800.00 &   & 300.00 & 4.00 & 125.80 & 429.80 & 1,229.80  \\
			5 & 0.00 & 120.00 & 1,480.00 & 1,600.00 &   & 300.00 & 4.00 & 206.80 & 510.80 & 2,110.80  \\
			6 & 0.00 & 120.00 & 1,810.00 & 1,930.00 &   & 60.00 & 1.00 & 219.00 & 280.00 & 2,210.00  \\
			7 & 0.00 & 200.00 & 850.00 & 1,050.00 &   & 120.00 & 27.00 & 162.80 & 309.80 & 1,359.80  \\
			8 & 1,200.00 & 200.00 & 320.00 & 1,720.00 &   & 480.00 & 33.00 & 23.60 & 536.60 & 2,256.60  \\
			9 & 0.00 & 200.00 & 1,940.00 & 2,140.00 &   & 120.00 & 33.00 & 202.00 & 355.00 & 2,495.00  \\
			10 & 1,200.00 & 200.00 & 2,470.00 & 3,870.00 &   & 480.00 & 27.00 & 336.40 & 843.40 & 4,713.40  \\
			11 & 0.00 & 200.00 & 2,470.00 & 2,670.00 &   & 120.00 & 27.00 & 336.40 & 483.40 & 3,153.40  \\
			12 & 1,200.00 & 200.00 & 1,940.00 & 3,340.00 &   & 480.00 & 33.00 & 202.00 & 715.00 & 4,055.00  \\
			\bottomrule                             
		\end{tabular}
	}
	\end{table}
	
	We also show the computational comparison of the model $\text{ACRP-S}^*$ with the model ACRP-A in Table \ref{Table:ArcModelCompare}. 
	The outcome demonstrates that model $\text{ACRP-S}^*$ performs better in terms of solution quality and overall run time.
	Although we get an integer gap equal to 0.00\% for all cases within model ACRP-A, the solution objectives of model $\text{ACRP-S}^*$ for all cases are equal to or better than those of model ACRP-A. Especially in large-scale scenarios with longer delay intervals (30 minutes), we observe bigger objective improvement for model $\text{ACRP-S}^*$. The reason is that delay interval of 30 minutes is not precise enough to get the optimal recovery solution in the arc-based model.
	
	\begin{table}[htbp]
		\centering\small
		\caption{\label{Table:ArcModelCompare} Comparison between Model ACRP-A and Model $\text{ACRP-S}^*$}
		\begin{tabular}{lrrrrrrrr}
			\toprule
			\multicolumn{1}{c}{} & \multicolumn{4}{c}{Model ACRP-A} & & \multicolumn{2}{c}{Model $\text{ACRP-S}^*$}  &                                                                                               \\
			\cmidrule{2-5} \cmidrule{7-8}
			\multirow{2}{*}{Scenario}       
			&\multicolumn{1}{c}{\begin{tabular}[c]{@{}c@{}}Max \\ delay\end{tabular}} &\multicolumn{1}{c}{\begin{tabular}[c]{@{}c@{}}Delay \\ interval\end{tabular}} 
			&\multicolumn{1}{c}{\begin{tabular}[c]{@{}c@{}}IP \\ Obj. \end{tabular}} 
			&\multicolumn{1}{c}{\begin{tabular}[c]{@{}c@{}}Total \\ time (s)\end{tabular}}
			&
			&\multicolumn{1}{c}{\begin{tabular}[c]{@{}c@{}}IP \\ Obj. \end{tabular}} 
			&\multicolumn{1}{c}{\begin{tabular}[c]{@{}c@{}}Total \\ time (s)\end{tabular}}
			&\multicolumn{1}{c}{\begin{tabular}[c]{@{}c@{}} Obj improve \\ of $\text{ACRP-S}^*$\end{tabular} } 
			\\
			\midrule
			1 & 240 & 5 & 632.20 & 803.0 &  & 632.20 & 3.4 & 0.00\%  \\
			2 & 240 & 5 & 738.40 & 3,610.0 &  & 738.40 & 2.9 & 0.00\%  \\
			3 & 240 & 5 & 1,604.60 & 753.0 &  & 1,604.60 & 3.5 & 0.00\%  \\
			4 & 240 & 5 & 1,229.80 & 3,076.0 &  & 1,229.80 & 3.9 & 0.00\%  \\
			5 & 240 & 5 & 2,110.80 & 790.0 &  & 2,110.80 & 3.0 & 0.00\%  \\
			6 & 240 & 5 & 2,210.00 & 2,642.0 &  & 2,210.00 & 3.1 & 0.00\%  \\
			7 & 180 & 30 & 1,754.40 & 4,831.0 &  & 1,359.80 & 727.4 & 22.49\%  \\
			8 & 180 & 30 & 2,499.20 & 8,906.0 &  & 2,256.60 & 314.8 & 9.71\%  \\
			9 & 180 & 30 & 2,805.20 & 3,320.0 &  & 2,495.00 & 1424.0 & 11.06\%  \\
			10 & 180 & 30 & 5,175.60 & 9,891.0 &  & 4,713.40 & 1131.1 & 8.93\%  \\
			11 & 180 & 30 & 3,615.60 & 5,005.0 &  & 3,153.40 & 918.6 & 12.78\%  \\
			12 & 180 & 30 & 4,365.20 & 7,909.0 &  & 4,055.00 & 665.5 & 7.11\%  \\			
			\bottomrule                             
		\end{tabular}
	\end{table}
	
	\subsection{Effectiveness of Valid Inequality Constraints}\label{Sec:ExperimentVI}
	
	As we stated in section \ref{Sec:ValidIneq}, including a set of valid inequality constraints into the string-based model ACRP-S would tighten boundaries and result in better recovery solution quality with a smaller integral gap.
	Inevitably, adding constraints for each cargo with each flight's departure time will increase the problem size, and the LP model's iterative solution time will increase as a result.
	
	\begin{table}[htbp]
		\centering\small
		\caption{\label{Table:ValidIneq}Comparison between Model ACRP-S and Model $\text{ACRP-S}^*$}
		\resizebox{\linewidth}{!}{  
		\begin{tabular}{lrrrrrrrrrrrrrr}
			\toprule
			\multicolumn{1}{c}{} & \multicolumn{6}{c}{Model ACRP-S Without Valid Inequality} & & \multicolumn{6}{c}{Model $\text{ACRP-S}^*$ With Valid Inequality}  &                                                        \\
			\cmidrule{2-7} \cmidrule{9-14}
			\multirow{2}{*}{Scenario}                      
			& \multicolumn{1}{c}{\begin{tabular}[c]{@{}c@{}}\# Col \\ Iter\end{tabular} } 
			& \multicolumn{1}{c}{\begin{tabular}[c]{@{}c@{}}LP \\ Obj.\end{tabular} } 
			& \multicolumn{1}{c}{\begin{tabular}[c]{@{}c@{}}LP \\ time\end{tabular} } 
			& \multicolumn{1}{c}{\begin{tabular}[c]{@{}c@{}}IP \\ Obj.\end{tabular}} 
			& \multicolumn{1}{c}{\begin{tabular}[c]{@{}c@{}}IP \\ time\end{tabular}} 
			& \multicolumn{1}{c}{\begin{tabular}[c]{@{}c@{}}LP\&IP \\ gap\end{tabular}} 
			&
			& \multicolumn{1}{c}{\begin{tabular}[c]{@{}c@{}}\# Col \\ Iter\end{tabular} } 
			& \multicolumn{1}{c}{\begin{tabular}[c]{@{}c@{}}LP \\ Obj.\end{tabular} } 
			& \multicolumn{1}{c}{\begin{tabular}[c]{@{}c@{}}LP \\ time\end{tabular} } 
			& \multicolumn{1}{c}{\begin{tabular}[c]{@{}c@{}}IP \\ Obj.\end{tabular}} 
			& \multicolumn{1}{c}{\begin{tabular}[c]{@{}c@{}}IP \\ time\end{tabular}} 
			& \multicolumn{1}{c}{\begin{tabular}[c]{@{}c@{}}LP\&IP \\ gap\end{tabular}} 
			& \multicolumn{1}{c}{\begin{tabular}[c]{@{}c@{}}LP \\ improve\end{tabular} }
			\\
			\midrule
			1 & 5 & 552.79 & 0.7 & 632.20 & 0.9 & 14.37\% &   & 6 & 594.86 & 2.8 & 632.20 & 3.4 & 6.28\% & 7.61\%  \\
			2 & 7 & 614.93 & 1.0 & 968.40 & 1.2 & 57.48\% &   & 9 & 706.33 & 2.7 & 738.40 & 2.9 & 4.54\% & 14.86\%  \\
			3 & 8 & 1,503.19 & 1.2 & 1,834.60 & 1.4 & 22.05\% &   & 8 & 1,572.59 & 3.2 & 1,604.60 & 3.5 & 2.04\% & 4.62\%  \\
			4 & 7 & 1,128.32 & 1.2 & 1,541.20 & 1.4 & 36.59\% &   & 9 & 1,192.39 & 3.7 & 1,229.80 & 3.9 & 3.14\% & 5.68\%  \\
			5 & 6 & 2,031.38 & 1.3 & 2,110.80 & 1.4 & 3.91\% &   & 7 & 2,073.45 & 2.8 & 2,110.80 & 3.0 & 1.80\% & 2.07\%  \\
			6 & 7 & 2,086.52 & 1.2 & 2,650.00 & 2.0 & 27.01\% &   & 8 & 2,177.92 & 2.8 & 2,210.00 & 3.1 & 1.47\% & 4.38\%  \\
			7 & 20 & 962.78 & 158.6 & 1,684.40 & 159.9 & 74.95\% &   & 37 & 1,310.21 & 719.5 & 1,359.80 & 727.4 & 3.78\% & 36.09\%  \\
			8 & 16 & 2,101.98 & 159.2 & 2,518.00 & 159.6 & 19.79\% &   & 14 & 2,221.94 & 311.0 & 2,256.60 & 314.8 & 1.56\% & 5.71\%  \\
			9 & 17 & 2,336.09 & 140.9 & 2,974.40 & 141.6 & 27.32\% &   & 25 & 2,455.42 & 1415.5 & 2,495.00 & 1424.0 & 1.61\% & 5.11\%  \\
			10 & 26 & 4,321.90 & 325.1 & 5,038.00 & 325.9 & 16.57\% &   & 33 & 4,658.83 & 1119.4 & 4,713.40 & 1131.1 & 1.17\% & 7.80\%  \\
			11 & 22 & 2,746.12 & 400.6 & 3,478.00 & 401.6 & 26.65\% &   & 32 & 3,098.87 & 908.9 & 3,153.40 & 918.6 & 1.76\% & 12.85\%  \\
			12 & 22 & 3,896.05 & 320.3 & 4,534.40 & 321.3 & 16.38\% &   & 27 & 4,015.39 & 660.5 & 4,055.00 & 665.5 & 0.99\% & 3.06\%  \\
			\bottomrule                             
		\end{tabular}
}	
\end{table}

In this section, we present the recovery solution result of model ACRP-S and model $\text{ACRP-S}^*$ in Table \ref{Table:ValidIneq}.
The solution approaches for the two models are identical as shown in Figure \ref{FigureFlowChart}.
The two models are identical and should produce the same optimal solution if we generate all potential aircraft strings and cargo itineraries. However, we only generate a subset of the total potential aircraft strings and cargo itineraries with the column-and-row generation approach. 
As a result, the relaxed LP solutions we obtained are lower bounds of the ideal optimal solution, whereas the IP solutions are upper bounds.
If we incorporate the valid inequality constraints into the model, the lower bounds of the LP solutions can be improved.

As we illustrated in Table \ref{Table:ValidIneq}, the addition of valid inequality constraints significantly improved the LP objective. The improvement ratio of LP objective is presented by "LP improve".
Additionally, we also obtained a small integer gap for model $\text{ACRP-S}^*$, which means that the LP solution is a good approximation for the ideal IP solution. If we solve model  ACRP-S by the column-and-row generation approach, we cannot guarantee that we will find enough "good" aircraft strings and cargo itineraries for the IP solution with a big integral gap. We present the details of recovery policies of the optimal IP solution for model ACRP-S in Table \ref{Table:NOVIdetail}. 
More cargo cancellations occurred in comparison to the results displayed in Table \ref{Table:IPdetail}, which is consistent with the illustration we provided in section \ref{Sec:ValidIneq}.

\begin{table}[htbp]
	\centering\small
	\caption{\label{Table:NOVIdetail}Recovery Policy Details of Model ACRP-S's IP Solutions}
	\begin{tabular}{lrrrrrrrr}
		\toprule
		\multicolumn{1}{c}{} & \multicolumn{3}{c}{Flight} & & \multicolumn{3}{c}{Cargo}                                                                                                                                  \\
		\cmidrule{2-4} \cmidrule{6-8}
		\multirow{2}{*}{Scenario}  &                    
		\multicolumn{1}{c}{\begin{tabular}[c]{@{}c@{}}\# flight \\ cancel\end{tabular} } &
		\multicolumn{1}{c}{\begin{tabular}[c]{@{}c@{}}\# aircraft \\ swap \end{tabular}} & 
		\multicolumn{1}{c}{\begin{tabular}[c]{@{}c@{}} Delay \\ minutes\end{tabular} } & 
		& 
		\multicolumn{1}{c}{\begin{tabular}[c]{@{}c@{}}\# cargo \\ cancel\end{tabular} } &
		\multicolumn{1}{c}{\begin{tabular}[c]{@{}c@{}}\# flight \\ change\end{tabular} } &
		\multicolumn{1}{c}{\begin{tabular}[c]{@{}c@{}} Delay \\ minutes\end{tabular} } 
		\\
		\midrule
		1 & 0 & 3 & 75 & 0 & 5 & 4 & 1,455  \\
		2 & 0 & 3 & 205 & 0 & 6 & 1 & 1,935  \\
		3 & 0 & 0 & 665 & 0 & 6 & 1 & 3,590  \\
		4 & 0 & 6 & 290 & 0 & 10 & 4 & 2,930  \\
		5 & 0 & 3 & 740 & 0 & 5 & 4 & 5,170  \\
		6 & 0 & 3 & 855 & 0 & 10 & 1 & 5,475  \\
		7 & 0 & 2 & 255 & 0 & 16 & 19 & 2,885  \\
		8 & 1 & 2 & 90 & 0 & 17 & 25 & 325  \\
		9 & 0 & 2 & 885 & 0 & 11 & 11 & 11,335  \\
		10 & 1 & 2 & 1,065 & 0 & 22 & 19 & 7,225  \\
		11 & 0 & 2 & 1,065 & 0 & 16 & 19 & 7,225  \\
		12 & 1 & 2 & 885 & 0 & 17 & 11 & 11,335  \\

		\bottomrule                             
	\end{tabular}
\end{table}

	We also noticed the increase of solution time for each iteration in model $\text{ACRP-S}^*$ in big cases. 
	Overall, the total solution time is a little bit longer for big scenarios with the combined disruption that occurred to schedule F284-A65 in practice. We will later show the results of the machine learning based column-and-row generation approach which is capable of reducing the problem size while also improving solution efficiency in section \cref{Sec:MLexperiments}.

	\subsection{Effectiveness of Integrated Model} \label{Sec:ExperimentSequencial}
	In order to demonstrate the effectiveness of integrated recovery, we compared model $\text{ACRP-S}^*$ with the sequential model, which solves flight recovery problems and cargo recovery problems sequentially. The sequential model, denoted by ACRP-SEQ, composes of two-stage models, ACRP-SEQ-I and ACRP-SEQ-II, as illustrated below: 

	\textbf{ACRP-SEQ-I: Flight Recovery}

	\begin{align}
		\min \ &\sum_{f \in F} c_{f} y_{f}+\sum_{a \in A} \sum_{l \in L} c_{a, l} x_{a, l} \label{SeqIObj}\\
		&\sum_{a \in A} \sum_{l \ni f} x_{a, l}+y_{f}=1&, &\forall f \in F \label{SeqICover}\\
		&\sum_{l \in L} x_{a, l} \leq 1&, &\forall a \in A \label{SeqIAft}\\
		&x_{a, l} \in\{0,1\}&, &\forall a \in A, \forall l \in L \label{SeqIX}\\
		&y_{f} \in\{0,1\}&, &\forall f \in F \label{SeqIY}
	\end{align}

	In model ACRP-SEQ-I, we only consider flight recovery decisions, ignoring cargo disruption issues. 
	The flights are assigned to the aircraft once the model ACRP-SEQ-I is solved, with each flight's cancellation and delay determined.
	As a result, the total capacity for each flight is also determined, too. 
	We denote the capacity for each flight $f$ as $Cap_{f}=\sum_{a \in A} \sum_{{l \owns f_t}} Cap_{a} x_{a, l}$ for simplification. And it is treated as a set of input parameters for the second stage. 

	The sequential model for the second stage is as follows:

\textbf{ACRP-SEQ-II: Cargo Recovery}
\begin{align}
	\min &\sum_{o \in O} c_{o} z_{o}+ \sum_{o \in O}\sum_{r \in R} c_{o,r} w_{o,r} \label{SeqIIObj}\\
	&Cap_f - \sum_{o\in O}\sum_{r\owns f}w_{o,r}\geq0,
	\ \forall f \in F
	\label{SeqIILegCap}
	\\
	& \sum_{r\in R}w_{o,r} + z_o = d_o, \forall o\in O \label{SeqIICargo}\\
	&z_{o}\in \mathbb{Z}, 0 \leq z_{o} \leq d_o, 	\ \forall o \in O
	\label{SeqIIZ}\\
	&w_{o,r} \in \mathbb{Z} 	\label{SeqIIW}
\end{align}

	The objective of the second stage model, ACRP-SEQ-II, is to re-assign cargoes to itineraries based on the rescheduled flights.
	The decision variables are the volume of cargo canceled and the volume of cargo transported by rescheduled cargo itineraries.
	
	We use the column generation approach to solve model ACRP-SEQ-I and model ACRP-SEQ-II sequentially. 
	To be more specific, in model ACRP-SEQ-I, we first include the original aircraft strings and then add better strings using the column generation approach.
	Better strings are generated for each aircraft by solving the multi-label shortest-path problem.
	When no more better air strings can be found, we solve the mixed-integer problem with an optimizing solver to get the integer solution of model ACRP-SEQ-I.
	With the integer solution of model ACRP-SEQ-I, we solve model ACRP-SEQ-II in a similar way.
	
	The result of model $\text{ACRP-S}^*$ in comparison with model ACRP-SEQ shows that, while the integrated model results in a solution with higher flight cost, the overall recovery cost is lower. The comparison details is shown in Table \ref{Table:SeqResult}. 
	In most scenarios, the integrated model $\text{ACRP-S}^*$ results in an objective improvement greater than 10\%.
	
	\begin{table}[htbp]
		\centering\small
		\caption{\label{Table:SeqResult}Comparison between Model ACRP-SEQ and Model $\text{ACRP-S}^*$}
		\begin{tabular}{lrrrrrrrr}
			\toprule
			\multicolumn{1}{c}{} & \multicolumn{3}{c}{Model ACRP-SEQ} & & \multicolumn{3}{c}{Model $\text{ACRP-S}^*$} &                                                          \\
			\cmidrule{2-4} \cmidrule{6-8}
			\multirow{2}{*}{Scenario}                      
			& \multicolumn{1}{c}{\begin{tabular}[c]{@{}c@{}}Flight \\ cost\end{tabular} } 
			& \multicolumn{1}{c}{\begin{tabular}[c]{@{}c@{}}Cargo \\ cost\end{tabular} } 
			& \multicolumn{1}{c}{\begin{tabular}[c]{@{}c@{}}Total \\ cost\end{tabular}} 
			&
			& \multicolumn{1}{c}{\begin{tabular}[c]{@{}c@{}}Flight \\ cost\end{tabular} } 
			& \multicolumn{1}{c}{\begin{tabular}[c]{@{}c@{}}Cargo \\ cost\end{tabular} }  
			& \multicolumn{1}{c}{\begin{tabular}[c]{@{}c@{}}Total \\ cost\end{tabular} } 
			&\multicolumn{1}{c}{\begin{tabular}[c]{@{}c@{}} Obj improve \\ of $\text{ACRP-S}^*$\end{tabular} } 
			\\
			\midrule
			
			1 & 60.00 & 1,139.20 & 1,199.20 &   & 270.00 & 362.20 & 632.20 & 47.28\%  \\
			2 & 500.00 & 678.40 & 1,178.40 &   & 600.00 & 138.40 & 738.40 & 37.34\%  \\
			3 & 1,330.00 & 504.60 & 1,834.60 &   & 1,400.00 & 204.60 & 1,604.60 & 12.54\%  \\
			4 & 560.00 & 1,446.80 & 2,006.80 &   & 800.00 & 429.80 & 1,229.80 & 38.72\%  \\
			5 & 1,390.00 & 1,287.80 & 2,677.80 &   & 1,600.00 & 510.80 & 2,110.80 & 21.17\%  \\
			6 & 1,830.00 & 820.00 & 2,650.00 &   & 1,930.00 & 280.00 & 2,210.00 & 16.60\%  \\
			7 & 590.00 & 1,094.40 & 1,684.40 &   & 1,050.00 & 309.80 & 1,359.80 & 19.27\%  \\
			8 & 1,430.00 & 1,310.80 & 2,740.80 &   & 1,720.00 & 536.60 & 2,256.60 & 17.67\%  \\
			9 & 1,850.00 & 1,124.40 & 2,974.40 &   & 2,140.00 & 355.00 & 2,495.00 & 16.12\%  \\
			10 & 3,410.00 & 1,628.00 & 5,038.00 &   & 3,870.00 & 843.40 & 4,713.40 & 6.44\%  \\
			11 & 2,210.00 & 1,268.00 & 3,478.00 &   & 2,670.00 & 483.40 & 3,153.40 & 9.33\%  \\
			12 & 3,050.00 & 1,484.40 & 4,534.40 &   & 3,340.00 & 715.00 & 4,055.00 & 10.57\%  \\
			
			\bottomrule                             
		\end{tabular}
	\end{table}

	\subsection{Effectiveness of Machine Learning Based Column-and-Row Generation Approach} \label{Sec:MLexperiments}
	In this section, we illustrate the effectiveness of the machine learning based column-and-row generation approach (denoted as \textit{ML-CRG}) with computational results.
	We first describe the details of our training data for machine learning, the feature selection process, as well as the prediction results. Then we present the computational results of \textit{ML-CRG} in comparison with \textit{CRG} as well as a heuristic approach.
	
	\subsubsection{Data and Features for Prediction Task} \label{SubSec:MLFeature}
	
	To build the decision tree model for the prediction task, we obtained 30 additional disruption datasets from the Chinese freighter airline mentioned above. Each data set contains a schedule with a 36-hour time-window and multiple fleets.
	The datasets include a variety of disruptions, such as previous flight delays, aircraft unavailability, airport closure, and so on. Then we solved the disruption scenarios with \textit{CRG} and recorded the solution details of flight connections during the column-and-row generation iterations. 
	
	Specifically, we record the LP result of each iteration for each flight connection during the solution process.
	In the LP solution, we can get a set of selected aircraft strings and a set of selected cargo itineraries with solution values greater than zero. 
	Each flight connection may be covered by several cargo itineraries with fractional LP solution value. Then we record all the departure times of the previous and successive flights in the selected itineraries, and check whether the connected flights make up a short through connection in any cargo itinerary.
	In the subsequent discussion, we designate the group of flight connections as the \textit{treatment group} if they are once recorded as short through connections through the column-and-row generation process. The remaining flight connections make up the \textit{control group}.
	With the 30 additional disruption datasets, we get 21,537 connections in total. The average amount of connections for each instance is 717.90, and the number varies between 463 and 1213. For all the connections, 905 are newly generated short through connections, almost 4.2\% of the total size. The average number of short through connections for each instance is 30.17, and the number varies between 14 and 45. The training data set is imbalanced with IR equals to 22.80.

	To predict the short through connection label $p_{con}$ for unobserved data set, it is critical to generate the features that are correlated with the observed label for the training data set $D.train$ while building the decision tree model.
	During each column-and-row iteration, we choose two critical features. The first feature is shortage of the connection time compared with standard cargo transshipment time (i.e., the connection time minus the standard cargo transshipment time), denoted as $h_1$.
	The second feature,  denoted as $h_2$, is the volume of uncovered cargo whose network contains the connection. 
	We propose that a flight connection is likely to become a short through connection with a shortage of connection time as well as some volume of cargo uncovered. 
	Furthermore, we observed that that connections that are solved as short through connections in one iteration are likely to be solved the same in following rounds. Therefore, throughout the feature creation process, we only collect the status of connections before becoming short through connections for the treatment group.
	
	Since we can observe a group of iteration outcomes for each connection, we aggregate the two features with minimum, maximum, mean, and standard variation, denoted as $h_1^{min}$, $h_1^{max}$, $h_1^{mean}$, $h_1^{sd}$ for the first feature, and $h_2^{min}$, $h_2^{max}$, $h_2^{mean}$, $h_2^{sd}$ for the second feature separately.
	To avoid multicollinearity among the features, we conduct a Spearman correlation test \citep{spearman1987}. The result shows that variable $h_1^{min}$, $h_1^{max}$ and $h_1^{mean}$ are highly correlated 
	(with Spearman correlation coefficient greater than 95\%). 
	We propose that greater shortage of the connection time are related to higher possibility of short through connections. Therefore, we keep feature $h_1^{min}$ and drop the other two features.
	
	\subsubsection{Predictive Result} \label{SubSec:MLResult}
\begin{figure}[htbp] 
		\centering 
		\includegraphics[width=0.6\textwidth]{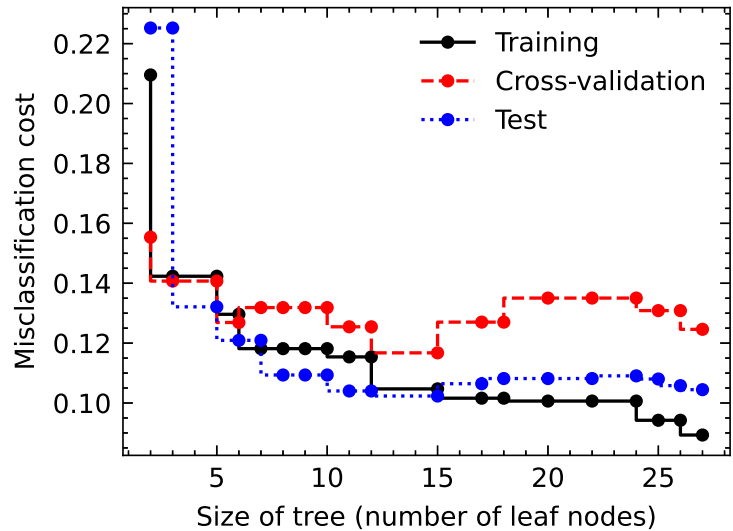} 
		\caption{The training, cross-validation, and test misclassification costs are shown as a function of the number of leaf nodes in the pruned tree.} 
		\label{FigureOverfit} 
	\end{figure}

	Now, we assess the prediction model's performance using the decision tree model. We first split the complete data set into a training set (80\%) and a testing set (20\%) using stratified sampling. After the entire tree was developed, we pruned the tree using the MCCP technique described above. Using a stratified 10-fold validation method, we obtain the average misclassification cost for the training and validation datasets.
	For every pruned subtree, we additionally computed the misclassification cost for $D.test$. in Figure \ref{FigureOverfit}, we display the training, cross-validation, and test misclassification costs are shown as a function of the number of leaf nodes in the pruned tree. The figure exhibits an over-fitting effect with the tree growing. That is, the misclassification cost for training decreased constantly as the tree size increased, but the cross-validation misclassification cost rose after the initial fall, with the minimum error at a tree size of 15. The cross-validation misclassification cost is a reasonable approximation of the test data set, because the test cost also dips down near the 15-node tree. Thus, we choose the pruned 15-node tree for the following integrated algorithm.
	
	Besides the misclassification cost, we are also interested in the sensitivity of the classification method. The sensitivity measures the portion of predicted short through connections among all actual short through connections. The sensitivity score with $D.test$ is 90.21\% in our case, which means we can capture more than 90\% percent of the total potential short through connections.
	
	
	\subsubsection{Computational Results of \textit{ML-CRG}} \label{SubSec:CGMLResult}
	
	We now describe the computational results for \textit{ML-CRG} in comparison with two competing solution methods. The first one is the column-and-row generation algorithm denoted as \textit{CRG}. The second one is a heuristic algorithm denoted as \textit{HEUR-CRG}. The heuristic approach works similar to \textit{ML-CRG}, while the only difference is that it takes the initial status of each flight connection as the prediction result of short through connections. The computational results for the scenarios in section \cref{Sec:ExperimentData} is shown in Table \ref{Table:MLCompare}.
	
	\begin{table}[htbp]
		\centering\small
		\caption{\label{Table:MLCompare}Comparison between algorithms CRG, HEUR-CRG and ML-CRG}
		\resizebox{\linewidth}{!}{  
			\begin{tabular}{lrrrrrrrrrrrrrrrrr}
				\toprule
				\multicolumn{1}{c}{} & \multicolumn{5}{c}{\textit{CRG}} & 
				\multicolumn{5}{c}{\textit{HEUR-CRG}} &  &   
				\multicolumn{5}{c}{\textit{ML-CRG}}                                                        \\
				\cmidrule{2-6} \cmidrule{8-12} \cmidrule{14-18}
				
				\multirow{2}{*}{Scenario}       
				& \multicolumn{1}{c}{\begin{tabular}[c]{@{}c@{}}\# SC \\ Const \end{tabular} } 
				& \multicolumn{1}{c}{\begin{tabular}[c]{@{}c@{}}LP \\ Obj.\end{tabular} } 
				& \multicolumn{1}{c}{\begin{tabular}[c]{@{}c@{}}LP \\ Time\end{tabular} } 
				& \multicolumn{1}{c}{\begin{tabular}[c]{@{}c@{}}IP \\ Obj.\end{tabular}} 
				& \multicolumn{1}{c}{\begin{tabular}[c]{@{}c@{}}IP \\ Time \end{tabular}} 
				&
				& \multicolumn{1}{c}{\begin{tabular}[c]{@{}c@{}}\# SC \\ Const \end{tabular} }             
				& \multicolumn{1}{c}{\begin{tabular}[c]{@{}c@{}}LP \\ Obj.\end{tabular} } 
				& \multicolumn{1}{c}{\begin{tabular}[c]{@{}c@{}}LP \\ Time\end{tabular} } 
				& \multicolumn{1}{c}{\begin{tabular}[c]{@{}c@{}}IP \\ Obj.\end{tabular}} 
				& \multicolumn{1}{c}{\begin{tabular}[c]{@{}c@{}}IP \\ Time \end{tabular}}
				&
				& \multicolumn{1}{c}{\begin{tabular}[c]{@{}c@{}}\# SC \\ Const \end{tabular} }            
				& \multicolumn{1}{c}{\begin{tabular}[c]{@{}c@{}}LP \\ Obj.\end{tabular} } 
				& \multicolumn{1}{c}{\begin{tabular}[c]{@{}c@{}}LP \\ Time\end{tabular} } 
				& \multicolumn{1}{c}{\begin{tabular}[c]{@{}c@{}}IP \\ Obj.\end{tabular}} 
				& \multicolumn{1}{c}{\begin{tabular}[c]{@{}c@{}}IP \\ Time \end{tabular}}
				\\
				\midrule
				1 & 9 & 594.86 & 2.8 & 632.20 & 3.4 &  & 9 & 682.19 & 1.7 & 682.20 & 1.9 &  & 10 & 594.86 & 2.8 & 632.20 & 3.0 \\ 
				2 & 16 & 706.33 & 2.7 & 738.40 & 2.9 &  & 10 & 824.32 & 1.0 & 824.40 & 1.1 &  & 15 & 706.32 & 1.7 & 738.40 & 1.9 \\ 
				3 & 16 & 1,572.59 & 3.2 & 1,604.60 & 3.5 &  & 9 & 1,604.59 & 1.5 & 1,604.60 & 1.7 &  & 16 & 1,572.59 & 3.1 & 1,604.60 & 3.4 \\ 
				4 & 12 & 1,192.39 & 3.7 & 1,229.80 & 3.9 &  & 10 & 1,359.72 & 2.6 & 1,359.80 & 2.8 &  & 12 & 1,192.39 & 3.3 & 1,229.80 & 3.6 \\ 
				5 & 18 & 2,073.45 & 2.8 & 2,110.80 & 3.0 &  & 9 & 2,160.79 & 1.2 & 2,160.80 & 1.3 &  & 16 & 2,073.45 & 1.8 & 2,110.80 & 2.3 \\ 
				6 & 23 & 2,177.92 & 2.8 & 2,210.00 & 3.1 &  & 9 & 2,295.92 & .9 & 2,296.00 & 1.0 &  & 21 & 2,177.92 & 2.0 & 2,210.00 & 2.3 \\ 
				7 & 898 & 1,310.21 & 719.5 & 1,359.80 & 727.4 &  & 297 & 1,430.06 & 141.7 & 1,450.80 & 145.7 &  & 844 & 1,310.22 & 487.3 & 1,359.80 & 493.7 \\ 
				8 & 683 & 2,221.94 & 311.0 & 2,256.60 & 314.8 &  & 191 & 2,334.97 & 99.6 & 2,347.60 & 101.8 &  & 544 & 2,221.94 & 183.3 & 2,256.60 & 186.7 \\ 
				9 & 1,054 & 2,455.42 & 1,415.5 & 2,495.00 & 1,424.0 &  & 324 & 2,573.38 & 250.8 & 2,586.00 & 257.2 &  & 916 & 2,455.42 & 897.8 & 2,495.00 & 904.4 \\ 
				10 & 1,015 & 4,658.83 & 1,119.4 & 4,713.40 & 1,131.1 &  & 279 & 4,784.94 & 188.8 & 4,804.40 & 193.2 &  & 928 & 4,658.84 & 940.1 & 4,713.40 & 951.9 \\ 
				11 & 954 & 3,098.87 & 908.9 & 3,153.40 & 918.6 &  & 308 & 3,224.98 & 170.8 & 3,244.40 & 176.2 &  & 824 & 3,098.87 & 702.2 & 3,153.40 & 712.4 \\ 
				12 & 1,076 & 4,015.39 & 660.5 & 4,055.00 & 665.5 &  & 274 & 4,133.35 & 102.1 & 4,146.00 & 104.4 &  & 650 & 4,015.38 & 451.9 & 4,055.00 & 456.1 \\ 
				
				\bottomrule           
			\end{tabular}
		}
	\end{table}
	
	The computational results show that algorithm \textit{CRG} is the slowest among the three algorithms. In Table \ref{Table:MLCompare}, "\# SC Const" is the short through connection constraints added to the master problem.
	The fastest algorithm, \textit{HEUR-CRG}, generates the fewest short through connection constraints, but results in solutions with greater recovery costs.
	The LP and IP objectives of \textit{ML-CRG} are equivalent to those of algorithm \textit{CRG}, indicating that that we added sufficient flight delay decisions of short through connections using the machine learning approach for the tested scenarios. 
	With machine learning, we save roughly 30\% of the time required to obtain the best LP solution for scenarios created from F284-A65 by adding fewer short through connections. 
	The majority of the big scenario solutions take less than 15 minutes, which is acceptable for recovery issues in real practice.
	As to Scenarios 1-6, the result shows little difference between \textit{CRG} and \textit{ML-CRG}, either in the number of generated short through connections or in the solution times. 
	
	\section{Conclusions}
	In this study, we introduced an integrated air cargo recovery problem that recovery flights, aircraft and cargoes simultaneously when disruptions happen.
	Two models based on the flight connection network are presented, one is the arc-based model, and the other is the string-based model. Different from the previous studies on the integrated problem, we not only consider aircraft and cargo re-routing decisions, but also make flight delay decisions for both aircraft and cargoes. We also consider short through connections when making flight delay decisions. By doing this, we utilize the benefits of short through connections and provide an integrated recovery solution with lower overall cost.
	
	
	
	The string-based model is solved via a machine learning based column-and-row generation approach. This approach uses the pricing strategy to limit the number of columns (aircraft strings and cargo itineraries) added into the model, thereby reducing the solution run times. 
	The pricing strategy has been shown to be effective in the literature, but due to the integer decision variables in the ACRP-S, there is a big gap between the linear relaxation solution and the mix-integer solution.
	Thus we added a set of valid inequality constraints to the string-based model in order to reduce the solution gap. 
	We also noticed that a large number of flight delay decision of short through connections are not selected in the column-and-row generation process. Therefore, to further improve the solution process, we developed a machine learning approach for short through connection prediction and integrate the prediction to the column-and-row generation process. This machine learning based column-and-row generation approach effectively reduced the number of flight delay decisions and constraints of short through connection that are added to the model, which speeds up the solution process and yields a good recovery solution in a reasonable amount of time.
	
	There are several directions for future work on our proposed model as well as solution approaches.
	One is to consider crew-recovery into the integrated model. Similar to cargo transshipment, it also costs a longer transit time for crew with sequential flights on different aircraft. This constraint can be formulated in a similar way as the short through connection constraint in our model. 
	Another possible extension is to recovery passenger aircraft, freighters and cargo itineraries simultaneously because many airlines use both passenger aircraft and freighters for cargo transportation.

	
	\ACKNOWLEDGMENT{%
		This study is supported by National Natural Science Foundation of China under Grant No.71825001.
	}
	
	%
	%
	%
	\begin{APPENDICES}
		\section{flight arc process in aircraft's network}
		\begin{breakablealgorithm}
			\caption{ Process flight arc $(i,j)$ in aircraft $a$'s network} 
			\label{Alg:Process Flight Arc} 
			\begin{algorithmic}[1] 
				\REQUIRE ~~\\ 
				Node $i$’s label set $B_{a,i}$, node $j$’s label set $B_{a,j}$, originally scheduled aircraft of flight $j$ and re-timing departure time set $T_j$. 
				\ENSURE ~~\\ 
				Updated label set $B_{a,j}$
				\FOR{ $ \forall \langle\bar{c}_{i}, t_{i}^{d}\rangle \in B_{a,i}$}
				\STATE Let $t_i$ be the re-timing departure time of flight $i$ with delay $t_{i}^{d}$
				\FOR {each departure time $t_{j} \in T_j$}
				\STATE $t_j=max\{t_j, t_i + fly\_time_i + turn\_time_{ij}\}$, $t_j^d = t_j - t_j^0$
				\IF{$t_{j}^{d}>maximum\_delay$
				}
				\STATE Continue loop
				\ENDIF
				\STATE Let $\bar{c}_{a,j} = \bar{c}_{a,i}+t_{j}^{d}\cdot delay\ cost-\alpha_{j} - Cap_{a} \gamma_{a,j}^{t_j} -\sum_{o\in O}\sum_{t'\in T_j\setminus t_j}\pi_{o,i,t'}$
				\IF {flight $j$ is not assigned to aircraft $a$ in the original plan}
				\STATE  $\bar{c}_{a,j} = \bar{c}_{a,j} +c_{a,i}^{swap}$
				\ENDIF
				\IF {flight arc $(i,j)$ with departure time $t_i$ and $t_j$ is a short through connection}
				\STATE  $\bar{c}_{a,j} = \bar{c}_{a,j} - Cap_a \eta_{(i^{t_i},j^{t_j})}$
				\ENDIF
				\IF {$\langle\bar{c}_{a,j}, t_{j}^{d}\rangle$ is not dominated by any label in $B_{a,j}$}
				\STATE  $B_{a,j} = B_{a,j} \cup \langle\bar{c}_{a,j}, t_{j}^{d}\rangle$ ; set $\langle\bar{c}_{a,j}, t_{j}^{d}\rangle$'s predecessor as $\langle\bar{c}_{a,i}, t_{i}^{d}\rangle$
				\STATE Delete labels dominated by $\langle\bar{c}_{a,j}, t_{j}^{d}\rangle$ in $B_{a,j}$
				\ENDIF
				\ENDFOR
				\ENDFOR
			\end{algorithmic}
		\end{breakablealgorithm}
		
		~\\
		\section{flight arc process in cargo's network}
		
		\begin{breakablealgorithm}
			\caption{ Process flight arc $(i,j)$ in cargo $o$'s network} 
			\label{Alg:Process Flight Arc_Cargo} 
			\begin{algorithmic}[1] 
				\REQUIRE ~~\\ 
				Node $i$’s label set $B_{o,i}$, node $j$’s label set $B_{o,j}$, original flight list scheduled to cargo $o$, and re-timing departure time set $T_j$
				\ENSURE ~~\\ 
				Updated label set $B_{o,j}$
				\FOR{ $\forall \langle\bar{c}_{o,i}, t_{i}^{d},\{sc\}_i\rangle \in B_{o,i}$}
				\STATE Let $t_i$ be the re-timing departure time of flight $i$ with delay $t_{i}^{d}$
				\STATE Let re-timing departure set $T_j^{new}$ of $j$  be $\{t_{i} + fly\_time_i+ turn\_time_{ij}, t_{i} + fly\_time_i+ trans\_time\}$
				\FOR {each departure time $t_{j} \in T_j \cup T_j^{new}$}
				\STATE $t_j^d = t_j - t_j^0$
				\IF {$t_j < t_i + fly\_time_i + turn\_time_{ij}$ or  $t_{j}^{d}>maximum\_delay$}
				\STATE Continue loop
				\ENDIF
				\STATE Let $\bar{c}_{o,j} = \bar{c}_{o,i}+ \gamma_{i,t} - \pi_{o,i,t}/d_o$,  $\{{sc}\}_j = \{{sc}\}_i$
				\IF {flight $j$ is not in the originally scheduled flights for cargo $o$}
				\STATE  $\bar{c}_{o,j} = \bar{c}_{o,j} +c_{o,i}^{change}$
				\ENDIF
				\IF {flight arc $(i,j)$ with departure time $t_i$ and $t_j$ is a short through connection}
				\STATE  $\bar{c}_{o,j} = \bar{c}_{o,j} +  \eta_{(i^{t_i},j^{t_j})}$,  $\{{sc}\}_j = \{{sc}\}_j \cup \{(i,j)\}$
				\ENDIF
				\IF {$\langle\bar{c}_{o,j}, t_{j}^{d},\{sc\}_j\rangle$ is not dominated by any label in $B_{o,j}$}
				\STATE  $B_{o,j} = B_{o,j} \cup \langle\bar{c}_{o,j}, t_{j}^{d},\{sc\}_j\rangle $ ; set $\langle\bar{c}_{o,j}, t_{j}^{d},\{sc\}_j\rangle$'s predecessor as $\langle\bar{c}_{o,i}, t_{i}^{d},\{sc\}_i\rangle$
				\STATE Delete labels dominated by $\langle\bar{c}_{o,j}, t_{j}^{d},\{sc\}_j\rangle$ in $B_{o,j}$
				\ELSE
				\STATE Delete label $\langle\bar{c}_{o,j}, t_{j}^{d},\{sc\}_j\rangle$ 
				\ENDIF
				\ENDFOR
				\ENDFOR
			\end{algorithmic}
		\end{breakablealgorithm}
		
	\end{APPENDICES}
	
	~\\
	
	
	\bibliographystyle{informs2014trsc} 
	\bibliography{mybibfile} 
	
	
\end{document}